\documentclass{elsart3-1-ppn}

\usepackage{amsmath,amsfonts,amssymb,amsxtra,setspace,xspace,graphicx,lmodern,fourier,graphicx}

\newtheorem{theorem}{Theorem}[section]
\newtheorem{lemma}[theorem]{Lemma}
\newtheorem{proposition}[theorem]{Proposition}
\newtheorem{corollary}[theorem]{Corollary}

\renewcommand{\theequation}{\arabic{equation}}
\newcommand{\C}{{\mathsf C}}
\newcommand{\R}{{\mathbb R}}
\renewcommand{\S}{{\mathbb S}}
\renewcommand{\L}{{\mathcal L_\alpha\,}}
\newcommand{\N}{{\mathbb N}}
\newcommand{\be}[1]{\begin{equation}\label{#1}}
\newcommand{\ee}{\end{equation}}
\renewcommand{\(}{\left(}
\renewcommand{\)}{\right)}
\newcommand{\ird}[1]{\int_{\R^d}{#1}\,dx}
\newcommand{\nrm}[2]{\left\|{#1}\right\|_{\mathrm L^{#2}(\R^d)}}
\newcommand{\irdmu}[1]{\int_{\R^d}#1\,d\mu}
\newcommand{\irdsph}[3]{\int_0^\infty\kern-5pt\int_{\S^{d-1}}#2\,{#1}^{\kern0.5pt #3-1}\,d{#1}\,d\omega}
\newcommand{\iring}[1]{\int_{ A_{(s,S)}}#1\,d\mu}
\newcommand{\icircle}[1]{\int_{\S^1}#1\,d\theta}
\newcommand{\idB}[2]{\int_{\partial B_{#1}}#2\,d\kern0.5pt\varsigma}
\newcommand{\isph}[1]{\int_{\S^{d-1}}#1\,d\omega}
\newcommand{\isphone}[2]{#1^{n-1}\(\int_{\S^{d-1}}#2\,d\kern0.5pt\omega\)(#1)}
\newcommand{\DD}[1]{\mathsf D_\alpha\kern0.5pt#1}
\newcommand{\dt}{d\kern-0.5pt t}

\newenvironment{proof}{\smallskip\par\noindent{\sl Proof.\/} }{\hfill$\square$}
\usepackage[english,francais]{babel}
\def\og{\leavevmode\raise.3ex\hbox{$\scriptscriptstyle\langle\!\langle$~}}
\def\fg{\leavevmode\raise.3ex\hbox{~$\!\scriptscriptstyle\,\rangle\!\rangle$}}

\usepackage{color}

\makeatletter\def\@textbottom{\vskip \z@ \@plus 20pt}\let\@texttop\relax\makeatother

\begin{document}
\centerline{}
\begin{frontmatter}\vspace*{-2cm}
\selectlanguage{english}
\title{Symmetry for extremal functions in subcritical Caffarelli-Kohn-Nirenberg inequalities}

\selectlanguage{english}
\author[Ceremade]{Jean Dolbeault},
\ead{dolbeaul@ceremade.dauphine.fr}
\author[Ceremade]{Maria J.~Esteban},
\ead{esteban@ceremade.dauphine.fr}
\author[ML]{Michael Loss}, \and
\ead{loss@math.gatech.edu}
\author[MM]{Matteo Muratori}
\ead{matteo.muratori@unipv.it}

\address[Ceremade]{Ceremade, UMR CNRS n$^\circ$~7534, Universit\'e Paris-Dauphine, PSL research university, Place de Lattre de Tassigny, 75775 Paris 16, France}
\address[ML]{School of Mathematics, Georgia Institute of Technology, Skiles Building, Atlanta GA 30332-0160, USA}
\address[MM]{Dipartimento di Matematica \emph{Felice Casorati}, Universit\`a degli Studi di Pavia, Via A.~Ferrata 5, 27100 Pavia, Italy}
\begin{abstract}
\selectlanguage{english}We use the formalism of the R\'enyi entropies to establish the symmetry range of extremal functions in a family of subcritical Caffarelli-Kohn-Nirenberg inequalities. By extremal functions we mean functions which realize the equality case in the inequalities, written with optimal constants. The method extends recent results on critical Caffarelli-Kohn-Nirenberg inequalities. Using heuristics given by a nonlinear diffusion equation, we give a variational proof of a symmetry result, by establishing a rigidity theorem: in the symmetry region, all positive critical points have radial symmetry and are therefore equal to the unique positive, radial critical point, up to scalings and multiplications. This result is sharp. The condition on the parameters is indeed complementary of the condition which determines the region in which symmetry breaking holds as a consequence of the linear instability of radial optimal functions. Compared to the critical case, the subcritical range requires new tools. The Fisher information has to be replaced by R\'enyi entropy powers, and since some invariances are lost, the estimates based on the Emden-Fowler transformation have to be modified.
\vskip 0.5\baselineskip

\selectlanguage{francais}
\vskip 0.5\baselineskip \noindent
{\bf Symm\'etrie des fonctions extr\'emales pour des in\'egalit\'es de Caffarelli-Kohn-Nirenberg sous-critiques}\par\noindent
Nous utilisons le formalisme des entropies de R\'enyi pour \'etablir le domaine de sym\'etrie des fonctions extr\'emales dans une famille d'in\'egalit\'es de Caffarelli-Kohn-Nirenberg sous-critiques. Par fonctions extr\'emales, il faut comprendre des fonctions qui r\'ealisent le cas d'\'egalit\'e dans les in\'egalit\'es \'ecrites avec des constantes optimales. La m\'ethode \'etend des r\'esultats r\'ecents sur les in\'egalit\'es de Caffarelli-Kohn-Nirenberg critiques. En utilisant une heuristique donn\'ee par une \'equation de diffusion non-lin\'eaire, nous donnons une preuve variationnelle d'un r\'esultat de sym\'etrie, gr\^ace \`a un th\'eor\`eme de rigidit\'e: dans la r\'egion de sym\'etrie, tous les points critiques positifs sont \`a sym\'etrie radiale et sont par cons\'equent \'egaux \`a l'unique point critique radial, positif, \`a une multiplication par une constante et \`a un changement d'\'echelle pr\`es. Ce r\'esultat est optimal. La condition sur les param\`etres est en effet compl\'ementaire de celle qui d\'efinit la r\'egion dans laquelle il y a brisure de sym\'etrie du fait de l'instabilit\'e lin\'eaire des fonctions radiales optimales. Compar\'e au cas critique, le domaine sous-critique n\'ecessite de nouveaux outils. L'information de Fisher doit \^etre remplac\'ee par l'entropie de R\'enyi, et comme certaines invariances sont perdues, les estimations bas\'ees sur la transformation d'Emden-Fowler doivent \^etre modifi\'ees.\medskip\noindent

\par\bigskip
\end{abstract}
\end{frontmatter}
\selectlanguage{english}

\section{A family of subcritical Caffarelli-Kohn-Nirenberg interpolation inequalities}\label{Sec:Intro}

With the norms
\[
\nrm w{q,\gamma}:=\(\ird{|w|^q\,|x|^{-\gamma}}\)^{1/q}\,,\quad\nrm wq:=\nrm w{q,0}\,,
\]
let us define $\mathrm L^{q,\gamma}(\R^d)$ as the space of all measurable functions $w$ such that $\nrm w{q,\gamma}$ is finite. Our functional framework is a space $\mathrm H^p_{\beta,\gamma}(\R^d)$ of functions $w\in\mathrm L^{p+1,\gamma}(\R^d)$ such that $\nabla w\in\mathrm L^{2,\beta}(\R^d)$, which is defined as the completion of the space $\mathcal D(\R^d\setminus\{0\})$ of the smooth functions on $\R^d$ with compact support in~$\R^d\setminus\{0\}$, with respect to the norm given by $\|w\|^2:=(p_\star-p)\,\nrm w{p+1,\gamma}^2+\nrm{\nabla w}{2,\beta}^2$.

Now consider the family of \emph{Caffarelli-Kohn-Nirenberg interpolation inequalities} given by\par\smallskip
\be{CKN}
\nrm w{2p,\gamma}\le\C_{\beta,\gamma,p}\,\nrm{\nabla w}{2,\beta}^\vartheta\,\nrm w{p+1,\gamma}^{1-\vartheta}\quad\forall\,w\in\mathrm H^p_{\beta,\gamma}(\R^d)\,.
\ee
\par\medskip\noindent Here the parameters $\beta$, $\gamma$ and $p$ are subject to the restrictions
\be{parameters}
d\ge2\,,\quad\gamma-2<\beta<\frac{d-2}d\,\gamma\,,\quad\gamma\in(-\infty,d)\,,\quad p\in\(1,p_\star\right]\quad\mbox{with}\quad p_\star:=\frac{d-\gamma}{d-\beta-2}
\ee
and the exponent $\vartheta$ is determined by the scaling invariance, \emph{i.e.},
\[\label{theta}
\vartheta=\frac{(d-\gamma)\,(p-1)}{p\,\big(d+\beta+2-2\,\gamma-p\,(d-\beta-2)\big)}\,.
\]
These inequalities have been introduced, among others, by L.~Caffarelli, R.~Kohn and L.~Nirenberg in~\cite{Caffarelli-Kohn-Nirenberg-84}. We observe that $\vartheta=1$ if $p=p_\star$, a case which has been dealt with in~\cite{DEL2015}, and we shall focus on the sub-critical case $p<p_\star$. Throughout this paper, $\C_{\beta,\gamma,p}$ denotes the optimal constant in \eqref{CKN}. We shall say that a function $w\in\mathrm H^p_{\beta,\gamma}(\R^d)$ is an \emph{extremal function} for~\eqref{CKN} if equality holds in the inequality.

\begin{figure}[ht]\begin{center}
\includegraphics[width=10cm]{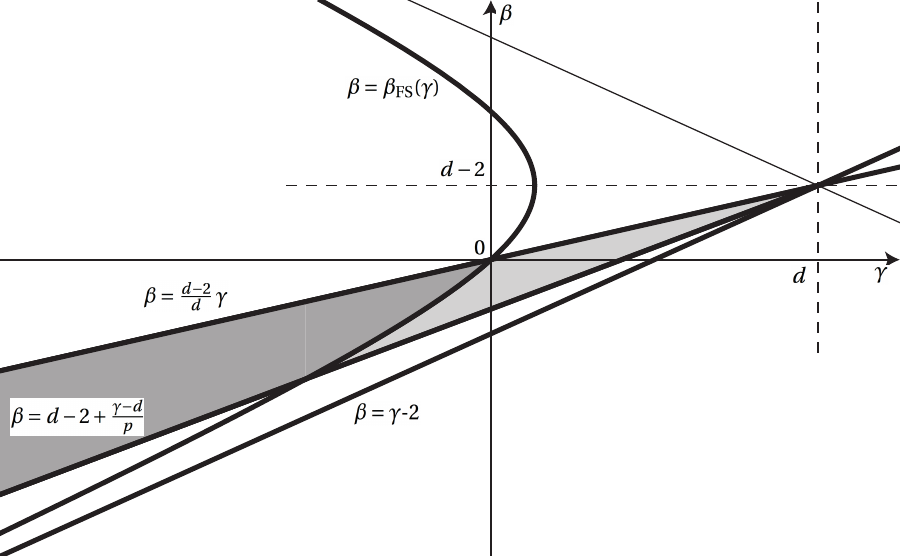}
\caption{\label{Fig}\emph{In dimension $d=4$, with $p=1.2$, the grey area corresponds to the cone determined by $d-2+(\gamma-d)/p\le\beta<(d-2)\,\gamma/d$ and $\gamma\in(-\infty,d)$ in~\eqref{parameters}. The light grey area is the region of symmetry, while the dark grey area is the region of symmetry breaking. The threshold is determined by the hyperbola $(d-\gamma)^2-(\beta-d+2)^2-4\,(d-1)=0$ or, equivalently $\beta=\beta_{\rm FS}(\gamma)$. Notice that the condition $p\le p_\star$ induces the restriction $\beta\ge d-2+(\gamma-d)/p$, so that the region of symmetry is bounded. The largest possible cone is achieved as $p\to1$ and is limited from below by the condition $\beta>\gamma-2$.}}
\end{center}\end{figure}

\medskip\emph{Symmetry} in~\eqref{CKN} means that the equality case is achieved by Aubin-Talenti type functions
\[
w_\star(x)=\(1+|x|^{2+\beta-\gamma}\)^{-1/(p-1)}\quad\forall\,x\in\R^d\,.
\]
On the contrary, there is \emph{symmetry breaking} if this is not the case, because the equality case is then achieved by a non-radial extremal function. It has been proved in~\cite{2016arXiv160208319B} that \emph{symmetry breaking} holds in~\eqref{CKN}~if
\be{set-symm-breaking}
\gamma<0\quad\mbox{and}\quad\beta_{\rm FS}(\gamma)<\beta<\frac{d-2}d\,\gamma
\ee
where
\[
\beta_{\rm FS}(\gamma):=d-2-\sqrt{(\gamma-d)^2-4\,(d-1)}\,.
\]
For completeness, we will give a short proof of this result in Section~\ref{Sec:SB}. Our main result shows that, under Condition~\eqref{parameters}, \emph{symmetry} holds in the complement of the set defined by~\eqref{set-symm-breaking}, which means that~\eqref{set-symm-breaking} is the sharp condition for \emph{symmetry breaking}. See Fig.~\ref{Fig}.
\par\smallskip\begin{theorem}\label{Thm:Main}{\sl Assume that~\eqref{parameters} holds and that
\be{Symmetry condition}
\beta\le\beta_{\rm FS}(\gamma)\quad\mbox{if}\quad\gamma<0\,.
\ee
Then the extremal functions for~\eqref{CKN} are radially symmetric and, up to a scaling and a multiplication by a constant, equal to $w_\star$.}\end{theorem}\par\smallskip
The result is slightly stronger than just characterizing the range of $(\beta,\gamma)$ for which equality in~\eqref{CKN} is achieved by radial functions. Actually our method of proof allows us to analyze the symmetry properties not only of extremal functions of~\eqref{CKN}, but also of all positive solutions in $\mathrm H^p_{\beta,\gamma}(\R^d)$ of the corresponding Euler-Lagrange equations, that is, up to a multiplication by a constant and a dilation, of
\be{ELeq}
-\,\mbox{div}\,\big(|x|^{-\beta}\,\nabla w\big)=|x|^{-\gamma}\,\big(w^{2p-1}-\,w^p\big)\quad\mbox{in} \quad \R^d \setminus \{ 0 \}\,.
\ee
\par\smallskip\begin{theorem}\label{Thm:Rigidity}{\sl Assume that~\eqref{parameters} and~\eqref{Symmetry condition} hold. Then all positive solutions to~\eqref{ELeq} in $\mathrm H^p_{\beta,\gamma}(\R^d)$ are radially symmetric and, up to a scaling and a multiplication by a constant, equal to $w_\star$.}\end{theorem}\par\smallskip
Up to a multiplication by a constant, we know that all non-trivial extremal functions for~\eqref{CKN} are non-negative solutions to~\eqref{ELeq}. Non-negative solutions to~\eqref{ELeq} are actually positive by the standard Strong Maximum principle. Theorem~\ref{Thm:Main} is therefore a consequence of Theorem~\ref{Thm:Rigidity}. In the particular case when $\beta=0$, the condition~\eqref{parameters} amounts to $d\ge2$, $\gamma\in (0,2)$, $p\in\big(1,(d-\gamma)/(d-2)\big]$, and~\eqref{CKN} can be written as
\[\label{CKN-beta0}
\nrm w{2p,\gamma}\le\C_{0,\gamma,p}\,\nrm{\nabla w}2^\vartheta\,\nrm w{p+1,\gamma}^{1-\vartheta}\quad\forall\,w\in\mathrm H^p_{0,\gamma}(\R^d)\,.
\]
In this case, we deduce from Theorem~\ref{Thm:Main} that symmetry always holds. This is consistent with a previous result ($\beta=0$ and $\gamma>0$, close to $0$) obtained in~\cite{DMN2015}. A few other cases were already known. The Caffarelli-Kohn-Nirenberg inequalities that were discussed in~\cite{DEL2015} correspond to the critical case $\theta=1$, $p=p_\star$ or, equivalently $\beta=d-2+(\gamma-d)/p$. Here by \emph{critical} we simply mean that $\nrm w{2p,\gamma}$ scales like $\nrm{\nabla w}{2,\beta}$. The limit case $\beta=\gamma-2$ and $p=1$, which is an endpoint for~\eqref{parameters}, corresponds to Hardy-type inequalities: there is no extremal function, but optimality is achieved among radial functions: see~\cite{0902}. The other endpoint is $\beta=(d-2)\,\gamma/d$, in which case $p_\star=d/(d-2)$. The results of Theorem~\ref{Thm:Main} also hold in that case with $p=p_\star=d/(d-2)$, up to existence issues: according to~\cite{Catrina-Wang-01}, either $\gamma\ge0$, symmetry holds and there exists a symmetric extremal function, or $\gamma<0$, and then symmetry is broken but there is no optimal function.

\medskip Inequality~\eqref{CKN} can be rewritten as an interpolation inequality with same weights on both sides using a change of variables. Here we follow the computations in~\cite{2016arXiv160208319B} (also see~\cite{DEL2015,dolbeault:hal-01286546}). Written in spherical coordinates for a function
\[
\widetilde w(r,\omega)=w(x)\,,\quad\mbox{with}\quad r=|x|\quad\mbox{and}\quad\omega=\frac x{|x|}\,,
\]
inequality~\eqref{CKN} becomes
\[
\(\irdsph r{|\widetilde w|^{2p}}{d-\gamma}\)^\frac1{2p}\le\C_{\beta,\gamma,p}\(\irdsph r{\left|\nabla \widetilde w\right|^2}{d-\beta}\)^\frac\vartheta2\(\irdsph r{|\widetilde w|^{p+1}}{d-\gamma}\)^\frac{1-\vartheta}{p+1}
\]
where $\left|\nabla\widetilde w\right|^2=\big|\tfrac{\partial\widetilde w}{\partial r}\big|^2+\tfrac1{r^2}\,\left|\nabla_{\kern-2pt\omega}\widetilde w\right|^2$ and $\nabla_{\kern-2pt\omega}\widetilde w$ denotes the gradient of $\widetilde w$ with respect to the angular variable $\omega\in\S^{d-1}$. Next we consider the change of variables $r\mapsto s=r^\alpha$,
\be{wv}
\widetilde w(r,\omega)=v(s,\omega)\quad\forall\,(r,\omega)\in\R^+\times\S^{d-1}
\ee
where $\alpha$ and $n$ are two parameters such that
\[
n=\frac{d-\beta-2}\alpha+2=\frac{d-\gamma}\alpha\,.
\]
Our inequality can therefore be rewritten as
\begin{multline*}
\(\irdsph s{|v|^{2p}}n\)^\frac1{2p}\\
\le\mathsf K_{\alpha,n,p}\(\irdsph s{\(\alpha^2\,\big|\tfrac{\partial v}{\partial s}\big|^2+\tfrac1{s^2}\,|\nabla_{\kern-2pt\omega} v|^2\)}n\)^\frac\vartheta2\(\irdsph s{|v|^{p+1}}n\)^\frac{1-\vartheta}{p+1}\,,
\end{multline*}
\[
\mbox{with}\quad\C_{\beta,\gamma,p}=\alpha^\zeta\,\mathsf K_{\alpha,n,p}\quad\mbox{and}\quad\zeta:=\frac\vartheta2+\frac{1-\vartheta}{p+1}-\frac1{2\,p}=\frac{(\beta+2-\gamma)\,(p-1)}{2\,p\,\big(d+\beta+2-2\,\gamma-p\,(d-\beta-2)\big)}\,.
\]
Using the notation
\[
\DD v=\(\alpha\,\frac{\partial v}{\partial s},\frac1s\,\nabla_{\kern-2pt\omega}v\)\,,
\]
with
\[\label{Eqn:alpha-n}
\alpha=1+\frac{\beta-\gamma}2\quad\mbox{and}\quad n=2\,\frac{d-\gamma}{\beta+2-\gamma}\,,
\]
Inequality~\eqref{CKN} is equivalent to a Gagliardo-Nirenberg type inequality corresponding to an artificial dimension~$n$ or, to be precise, to a Caffarelli-Kohn-Nirenberg inequality with weight $|x|^{n-d}$ in all terms. Notice that
\[
p_\star=\frac n{n-2}\,.
\]
\par\smallskip\begin{corollary}\label{Cor:Main}{\sl Assume that $\alpha$, $n$ and $p$ are such that
\[\label{parameters1}
d\ge2\,,\quad\alpha>0\,,\quad n>d\quad\mbox{and}\quad p\in\(1,p_\star\right]\,.
\]
Then the inequality
\be{CKN1}
\nrm v{2p,d-n}\le\mathsf K_{\alpha,n,p}\,\nrm{\DD v}{2,d-n}^\vartheta\,\nrm v{p+1,d-n}^{1-\vartheta}\quad\forall\,v\in\mathrm H^p_{d-n,d-n}(\R^d)\,,
\ee
holds with optimal constant $\mathsf K_{\alpha,n,p}=\alpha^{-\kern0.5pt\zeta}\,\C_{\beta,\gamma,p}$ as above and optimality is achieved among radial functions if and only if
\be{SymmetryCondition1}
\alpha\le\alpha_{\rm FS}\quad\mbox{with}\quad\alpha_{\rm FS}:=\sqrt{\frac{d-1}{n-1}}\,.
\ee
When symmetry holds, optimal functions are equal, up to a scaling and a multiplication by a constant, to
\[
v_\star(x):=\(1+|x|^2\)^{-1/(p-1)}\quad\forall\,x\in\R^d\,.
\]
}\end{corollary}\par\smallskip
We may notice that neither $\alpha_{\rm FS}$ nor $\beta_{\rm FS}$ depend on $p$ and that the curve $\alpha=\alpha_{\rm FS}$ determines the same threshold for the symmetry breaking region as in the critical case $p=p_\star$. In the case $p=p_\star$, this curve was found by V.~Felli and M.~Schneider, who proved in~\cite{Felli-Schneider-03} the linear instability of all radial critical points if $\alpha>\alpha_{\rm FS}$. When $p=p_\star$, symmetry holds under Condition~\eqref{SymmetryCondition1} as was proved in~\cite{DEL2015}. Our goal is to extend this last result to the subcritical regime $p\in(1,p_\star)$.

\medskip The change of variables $s=r^\alpha$ is an important intermediate step, because it allows to recast the problem as a more standard interpolation inequality in which the \emph{dimension} $n$ is, however, not necessarily an integer. Actually $n$ plays the role of a dimension in view of the scaling properties of the inequalities and, with respect to this \emph{dimension}, they are critical if $p=p_\star$ and sub-critical otherwise. The critical case $p=p_\star$ has been studied in~\cite{DEL2015} using tools of entropy methods, a critical fast diffusion flow and, in particular, a reformulation in terms of a \emph{generalized Fisher information}. In the subcritical range, we shall replace the entropy by a \emph{R\'enyi entropy power} as in~\cite{MR3200617,1501}, and make use of the corresponding fast diffusion flow. As in~\cite{DEL2015}, the flow is used only at heuristic level in order to produce a well-adapted test function. The core of the method is based on the Bakry-Emery computation, also known as the \emph{carr\'e du champ method}, which is well adapted to optimal interpolation inequalities: see for instance~\cite{MR3155209} for a general exposition of the method and~\cite{MR3229793,dolbeault:hal-01206975} for its use in presence of nonlinear flows. Also see~\cite{MR1853037} for earlier considerations on the Bakry-Emery method applied to nonlinear flows and related functional inequalities in unbounded domains. However, in non-compact manifolds and in presence of weights, integrations by parts have to be justified. In the critical case, one can rely on an additional invariance to use an Emden-Fowler transformation and rewrite the problem as an autonomous equation on a cylinder, which simplifies the estimates a lot. In the subcritical regime, estimates have to be adapted since after the Emden-Fowler transformation,
 the problem in the cylinder is no longer autonomous. 
 
This paper is organized as follows. We recall the computations which characterize the linear instability of radially symmetric minimizers in Section~\ref{Sec:SB}. In Section~\ref{RenyiStrategy}, we expose the strategy for proving symmetry in the subcritical regime when there are no weights. Section~\ref{Sec:BE} is devoted to the Bakry-Emery computation applied to R\'enyi entropy powers, in presence of weights. This provides a proof of our main results, if we admit that no boundary term appears in the integrations by parts in Section~\ref{Sec:BE}. To prove this last result, regularity and decay estimates of positive solutions to~\eqref{ELeq} are established in Section~\ref{Sec:RegDecay}, which indeed show that no boundary term has to be taken into account (see Proposition~\ref{Prop:b}).

\section{Symmetry breaking}\label{Sec:SB}

For completeness, we summarize known results on symmetry breaking for~\eqref{CKN}. Details can be found in~\cite{2016arXiv160208319B}. With the notations of Corollary~\ref{Cor:Main}, let us define the functional
\[
\mathcal J[v]:=\vartheta\,\log\(\nrm{\DD v}{2,d-n}\)+(1-\vartheta)\,\log\(\nrm v{p+1,d-n}\)+\log\mathsf K_{\alpha,n,p}-\log\(\nrm v{2p,d-n}\)
\]
obtained by taking the difference of the logarithm of the two terms in~\eqref{CKN1}. Let us define $d\mu_\delta:=\mu_\delta(x)\,dx$, where
\[
\mu_\delta(x):=\frac1{(1+|x|^2)^\delta}\,.
\]
Since~$v_\star$ as defined in Corollary~\ref{Cor:Main} is a critical point of $\mathcal J$, a Taylor expansion at order $\varepsilon^2$ shows that
\[
\nrm{\DD v_\star}{2,d-n}^2\,\mathcal J\big[v_\star+\varepsilon\,\mu_{\delta/2}\,f\big]=\tfrac12\,\varepsilon^2\,\vartheta\,\mathcal Q[f]+o(\varepsilon^2)
\]
with $\delta=\frac{2\,p}{p-1}$ and
\[
\mathcal Q[f]=\irdmu{|\DD f|^2\,|x|^{n-d}}_\delta-\,\frac{4\,p\,\alpha^2}{p-1}\irdmu{|f|^2\,|x|^{n-d}}_{\delta+1}\,.
\]
The following \emph{Hardy-Poincar\'e inequality} has been established in~\cite{2016arXiv160208319B}.
\par\smallskip\begin{proposition}\label{Prop:SpectralGap}{\sl Let $d\ge2$, $\alpha\in(0,+\infty)$, $n>d$ and $\delta\ge n$. Then
\be{Hardy-Poincare1}
\irdmu{|\DD f|^2\,|x|^{n-d}}_\delta\ge\Lambda\irdmu{|f|^2\,|x|^{n-d}}_{\delta+1}
\ee
holds for any $f\in\mathrm L^2(\R^d,|x|^{n-d}\,d\mu_{\delta+1})$, with $ \DD f \in \mathrm L^2(\R^d,|x|^{n-d}\,d\mu_{\delta}) $, such that $\irdmu{f\,|x|^{n-d}}_{\delta+1}=0$, with an optimal constant $\Lambda$ given by
\[\label{Eqn:SpectralGap}
\Lambda=\left\{\begin{array}{rl}
2\,\alpha^2\,(2\,\delta-n)\quad&\mbox{if}\quad0<\alpha^2\le\frac{(d-1)\,\delta^2}{n\,(2\,\delta-n)\,(\delta-1)}\,,\\[6pt]
2\,\alpha^2\,\delta\,\eta\quad&\mbox{if}\quad\alpha^2>\frac{(d-1)\,\delta^2}{n\,(2\,\delta-n)\,(\delta-1)}\,,
\end{array}
\right.
\]
where $\eta$ is the unique positive solution to
\[
\eta\,(\eta+n-2)=\frac{d-1}{\alpha^2}\,.
\]
Moreover, $\Lambda$ is achieved by a non-trivial eigenfunction corresponding to the equality in~\eqref{Hardy-Poincare1}. If $\alpha^2>\frac{(d-1)\,\delta^2}{n\,(2\,\delta-n)\,(\delta-1)}$, the eigenspace is generated by $\varphi_i(s,\omega)=s^\eta\,\omega_i$, with $i=1$, $2$,\ldots$d$ and the eigenfunctions are not radially symmetric, while in the other case the eigenspace is generated by the radially symmetric eigenfunction $\varphi_0(s,\omega)=s^2-\frac n{2\,\delta-n}$.}\end{proposition}\par\smallskip
As a consequence, $\mathcal Q$ is a nonnegative quadratic form if and only if $\frac{4\,p\,\alpha^2}{p-1}\le\Lambda$. Otherwise, $\mathcal Q$ takes negative values, and a careful analysis shows that symmetry breaking occurs in~\eqref{CKN} if
\[
2\,\alpha^2\,\delta\,\eta<\frac{4\,p\,\alpha^2}{p-1}\quad\Longleftrightarrow\quad\eta<1\,,
\]
which means
\[
\frac{d-1}{\alpha^2}=\eta\,(\eta+n-2)<n-1\,,
\]
and this is equivalent to $\alpha>\alpha_{\rm FS}$.

\section{The strategy for proving symmetry without weights}\label{RenyiStrategy}

Before going into the details of the proof we explain the strategy for the case of the Gagliardo-Nirenberg inequalities without weights. There are several ways to compute the optimizers, and the relevant papers are~\cite{MR1940370,MR1777035,MR1986060,MR1853037,MR3155209,1501} (also see additional references therein). The inequality is of the form
\be{ordinary}
\nrm w{2p}\le\C_{0,0,p}\,\nrm{\nabla w}2^\vartheta\,\nrm w{p+1}^{1-\vartheta}\quad\mbox{with}\quad1<p<\frac d{d-2}
\ee
and
\[
\vartheta=\frac{d\,(p-1)}{p\,\big(d+2-p\,(d-2)\big)}\,.
\]

It is known through the work in~\cite{MR1940370} that the optimizers of this inequality are, up to multiplications by a constant, scalings and translations, given by
\[
w_\star(x)=\(1+|x|^2\)^{-\frac1{p-1}}\quad\forall\,x\in\R^d\,.
\]
In our perspective, the idea is to use a version of the \emph{carr\'e du champ} or \emph{Bakry-Emery method} introduced in~\cite{Bakry-Emery85}: by differentiating a relevant quantity along the flow, we recover the inequality in a form which turns out to be sharp. The version of the \emph{carr\'e du champ} we shall use is based on the \emph{R\'enyi entropy powers} whose concavity as a function of $t$ has been studied by M.~Costa in~\cite{MR823597} in the case of linear diffusions (see~\cite{MR3200617} and references therein for more recent papers). In~\cite{MR1768665}, C.~Villani observed that the \emph{carr\'e du champ} method gives a proof of the logarithmic Sobolev inequality in the Blachman-Stam form, also known as the Weissler form: see~\cite{MR0188004,MR479373}. G.~Savar\'e and G.~Toscani observed in~\cite{MR3200617} that the concavity also holds in the nonlinear case, which has been used in~\cite{1501} to give an alternative proof of the Gagliardo-Nirenberg inequalities, that we are now going to sketch. 

The first step consists in reformulating the inequality in new variables. We set
\[
u=w^{2p}\,,
\]
which is equivalent to $w=u^{m-1/2}$, and consider the flow given by
\be{ordinaryflow}
\frac{\partial u}{\partial t}=\Delta u^m\,,
\ee
where $m$ is related to $p$ by
\[
p=\frac1{2\,m-1}\,.
\]
The inequalities $1<p<\frac d{d-2}$ imply that
\be{ordinaryrange}
1-\frac1d<m<1\,.
\ee
For some positive constant $\kappa>0$, one easily finds that the so-called Barenblatt-Pattle functions
\[\label{ordinaryselfsimilar}
u_\star(t,x)=\kappa^d\,t^{-\frac d{d\,m-d+2}}\,w_\star^{2p}\(\kappa\,t^{-\frac1{d\,m-d+2}}\,x\)=\(a+b\,|x|^2\)^{-\frac1{1-m}}
\]
are self-similar solutions of \eqref{ordinaryflow}, where $a=a(t)$ and $b=b(t)$ are explicit. Thus, we see that $w_\star=u_\star^{m-1/2}$ is an optimizer for~\eqref{ordinary} for all $t$ and it makes sense to rewrite~\eqref{ordinary} in terms of the function $u$. Straightforward computations show that~\eqref{ordinary} can be brought into the form
\be{ordinaryu}
\(\ird u\)^{(\sigma+1)\,m-1}\le C\,\mathcal E^{\sigma-1}\,\mathcal I\quad\mbox{where}\quad\sigma=\frac2{d\,(1-m)}-1
\ee
for some constant $C$ which does not depend on $u$, where
\[
\mathcal E:=\ird{u^m}
\]
is a \emph{generalized Ralston-Newman entropy}, also known in the literature as \emph{Tsallis entropy}, and
\[
\mathcal I:=\ird{u\,|\nabla\mathsf P|^2}
\]
is the corresponding \emph{generalized Fisher information}. Here we have introduced the \emph{pressure variable}
\[
\mathsf P=\frac m{1-m}\,u^{m-1}\,.
\]
The \emph{R\'enyi entropy power} is defined by
\[
\mathcal F:=\mathcal E^\sigma
\]
as in~\cite{MR3200617,1501}. With the above choice of $\sigma$, $\mathcal F$ is an affine function of $t$ if $u=u_\star$. For an arbitrary solution of~\eqref{ordinaryflow}, we aim at proving that it is a concave function of $t$ and that it is affine if and only if $u=u_\star$. For further references on related issues see~\cite{MR1940370,MR2282669}. Note that one of the motivations for choosing the variable $\mathsf P$ is that it has a particular simple form for the self-similar solutions, namely
\[
\mathsf P_{\kern -1pt\star}=\frac m{1-m}\(a+b\,|x|^2\)\,.
\]

Differentiating $\mathcal E$ along the flow~\eqref{ordinaryflow} yields
\[
\mathcal E'=(1-m)\,\mathcal I\,,
\]
so that
\[
\mathcal F'=\sigma\,(1-m)\,\mathcal G\quad\mbox{with}\quad\mathcal G:=\mathcal E^{\sigma-1}\,\mathcal I\,.
\]
More complicated is the derivative for the Fisher information:
\[
\mathcal I'=-\,2\ird{u^m\,\left[{\rm Tr}\(\(\mathrm{Hess}\,\mathsf P-\tfrac1d\,\Delta\mathsf P\,{\rm Id}\)^2\)+\(m-1+\tfrac1d\)(\Delta\mathsf P)^2\right]}\,.
\]
Here $\mathrm{Hess}\,\mathsf P$ and $\mathrm{Id}$ are respectively the Hessian of $\mathsf P$ and the $(d\times d)$ identity matrix. The computation can be found in~\cite{1501}. Next we compute the second derivative of the \emph{R\'enyi entropy power} $\mathcal F$ with respect to $t$:
\[
\frac{(\mathcal F)''}{\sigma\,\mathcal E^\sigma}=(\sigma-1)\,\frac{\mathcal E'^2}{\mathcal E^2}+\frac{\mathcal E''}{\mathcal E}=(\sigma-1)\,(1-m)^2\,\frac{\mathcal I^2}{\mathcal E^2}+(1-m)\,\frac{\mathcal I'}{\mathcal E}=:(1-m)\,\mathcal H\,.
\]
With $\sigma=\frac2d\,\frac1{1-m}-1$, we obtain
\be{ordinaryexpression}
\mathcal H=-\,2\,\left\langle{\rm Tr}\(\(\mathrm{Hess}\,\mathsf P-\tfrac1d\,\Delta\mathsf P\,{\rm Id}\)^2\)\right\rangle+(1-m)\,(1-\sigma)\,\left\langle\(\Delta\mathsf P-\langle\Delta\mathsf P\rangle\)^2\right\rangle
\ee
where we have used the notation
\[
\langle A\rangle:=\frac{\ird{u^m\,A}}{\ird{u^m}}\,.
\]
Note that by~\eqref{ordinaryrange}, we have that $\sigma >1$ and hence we find that $\mathcal F''=(\mathcal E^\sigma)''\le0$, which also means that $\mathcal G=\mathcal E^{\sigma-1}\,\mathcal I$ is a non-increasing function. In fact it is strictly decreasing unless $\mathsf P$ is a polynomial function of order two in~$x$ and it is easy to see that the expression~\eqref{ordinaryexpression} vanishes precisely when $\mathsf P$ is of the form $a+b\,|x-x_0|^2$, where $a$, $b\in\R$, $x_0\in\R^d$ are constants (but $a$ and $b$ may still depend on $t$).

Thus, while the left side of~\eqref{ordinaryu} stays constant along the flow, the right side decreases. In~\cite{1501} it was shown that the right side decreases towards the value given by the self-similar solutions~$u_\star$ and hence proves~\eqref{ordinary} in the sharp form. In our work we pursue a different tactic. The variational equation for the optimizers of~\eqref{ordinary} is given~by
\[
-\Delta w=a\,w^{2\,p-1}-b\,w^p\,.
\]
A straightforward computation shows that this can be written in the form
\[
2\,m\,u^{m-2}\,{\rm div}\big(u\,\nabla\mathsf P\big)+|\nabla\mathsf P|^2+c_1\,u^{m-1}=c_2
\]
for some constants $c_1$, $c_2$ whose precise values are explicit. This equation can also be interpreted as the variational equation for the sharp constant in~\eqref{ordinaryu}. Hence, multiplying the above equation by $\Delta u^m$ and integrating yields
\[
\ird{\left[2\,m\,u^{m-2}\,{\rm div}\big(u\,\nabla\mathsf P\big)+|\nabla\mathsf P|^2\right]\,\Delta u^m}+c_1\ird{u^{m-1}\,\Delta u^m}=c_2\ird{\Delta u^m}=0\,.
\]
We recover the fact that, in the flow picture, $\mathcal H$ is, up to a positive factor, the derivative of $\mathcal G$ and hence vanishes. From the observations made above we conclude that $\mathsf P$ must be a polynomial function of order two in~$x$. In this fashion one obtains more than just the optimizers, namely a classification of all positive solutions of the variational equation. The main technical problem with this method is the justification of the integrations by parts, which in the case at hand, without any weight, does not offer great difficulties: see for instance~\cite{MR1853037}. This strategy can also be used to treat the problem with weights, which will be explained next. Dealing with weights, however, requires some special care as we shall see.

\section{The Bakry-Emery computation and R\'enyi entropy powers in the weighted case}\label{Sec:BE}

Let us adapt the above strategy to the case where there are weights in all integrals entering into the inequality, that is, let us deal with inequality~\eqref{CKN1} instead of inequality~\eqref{ordinary}. In order to define a new, well-adapted fast diffusion flow, we introduce the diffusion operator $\L=-\,\mathsf D_\alpha^*\,\DD$, which is given in spherical coordinates~by
\[
\L u=\alpha^2\(u''+\frac{n-1}s\,u'\)+\frac1{s^2}\,\Delta_\omega\,u
\]
where $\Delta_\omega$ denotes the Laplace-Betrami operator acting on the $(d-1)$-dimensional sphere $\S^{d-1}$ of the angular variables, and $'$ denotes here the derivative with respect to $s$. Consider the fast diffusion equation
\be{FDE}
\frac{\partial u}{\partial t}=\L u^m
\ee
in the subcritical range $1-\frac1n<m=1-\frac1\nu<1$. The exponents $m$ in~\eqref{FDE} and $p$ in~\eqref{CKN1} are related as in Section~\ref{RenyiStrategy} by
\[
p=\frac1{2\,m-1}\quad\Longleftrightarrow\quad m=\frac{p+1}{2\,p}
\]
and $\nu$ is defined by
\[
\nu:=\frac1{1-m}\,.
\]
We consider the \emph{Fisher information} defined as
\[
\mathcal I[\mathsf P]:=\irdmu{u\,|\DD\mathsf P|^2}\quad\mbox{with}\quad\mathsf P=\frac m{1-m}\,u^{m-1}\quad\mbox{and}\quad d\mu=s^{n-1}\,ds\,d\omega=s^{n-d}\,dx\,.
\]
Here $\mathsf P$ is the \emph{pressure variable}. Our goal is to prove that $\mathsf P$ takes the form $a+b\,s^2$, as in Section~\ref{RenyiStrategy}. It is useful to observe that~\eqref{FDE} can be rewritten as
\[
\frac{\partial u}{\partial t}=\mathsf D_\alpha^*\(u\,\DD\mathsf P\)
\]
and, in order to compute $\frac{d\mathcal I}{\dt}$, we will also use the fact that $\mathsf P$ solves
\be{Eqn:p2}
\frac{\partial\mathsf P}{\partial t}=(1-m)\,\mathsf P\,\L\mathsf P-\,|\DD\mathsf P|^2\,.
\ee

\subsection{First step: computation of~$\frac{d\kern-0.5pt\mathcal I}{\dt}$} Let us define
\[
\mathcal K[\mathsf P]:=\mathcal A[\mathsf P]-(1-m)\,(\L\mathsf P)^2\quad\mbox{where}\quad\mathcal A[\mathsf P]:=\frac12\,\L\,|\DD\mathsf P|^2-\,\DD\mathsf P\cdot\DD\L\mathsf P
\]
and, on the boundary of the centered ball $B_s$ of radius $s$, the boundary term 
\begin{multline}\label{b}
\mathsf b(s):=\idB s{\Big(\tfrac\partial{\partial s}\(\mathsf P^{\frac m{m-1}}\,|\DD\mathsf P|^2\)-\,2\,(1-m)\,\mathsf P^{\frac m{m-1}}\,\mathsf P'\,\L\mathsf P\Big)}\\
=\isphone s{\Big(\tfrac\partial{\partial s}\(\mathsf P^{\frac m{m-1}}\,|\DD\mathsf P|^2\)-\,2\,(1-m)\,\mathsf P^\frac m{m-1}\,\mathsf P'\,\L\mathsf P\Big)}\,,
\end{multline}
where by $ d\varsigma = s^{n-1}\,d\omega $ we denote the standard Hausdorff measure on $ \partial B_s $.
\par\smallskip\begin{lemma}\label{Lem:DerivFisherL}{\sl If $u$ solves~\eqref{FDE} and if 
\be{FisherDerBC}
\lim_{s\to0_+}\mathsf b(s)=\lim_{S\to+\infty}\mathsf b(S)=0\,,
\ee
then,
\be{BLWL}
\frac d{\dt}\,\mathcal I[\mathsf P]=-\,2\irdmu{\mathcal K[\mathsf P]\,u^m}\,.
\ee
}\end{lemma}\par\smallskip
\begin{proof} For $0<s<S<+\infty$, let us consider the set $A_{(s,S)}:=\left\{x\in\R^d\,:\,s<|x|<S\right\}$, so that $\partial A_{(s,S)}=\partial B_s\cup\partial B_S$. Using~\eqref{FDE} and~\eqref{Eqn:p2}, we can compute
\begin{eqnarray*}
\hspace*{2cm}&&\hspace*{-2cm}\frac d{\dt}\iring{u\,|\DD\mathsf P|^2}\\
&=&\iring{\frac{\partial u}{\partial t}\,|\DD\mathsf P|^2}+\,2\iring{u\,\DD\mathsf P\cdot\DD\frac{\partial\mathsf P}{\partial t}}\\
&=&\iring{\L(u^m)\,|\DD\mathsf P|^2}+\,2\iring{u\,\DD\mathsf P\cdot\DD\Big((1-m)\,\mathsf P\,\L\mathsf P-|\DD\mathsf P|^2\Big)}\\
&=&\iring{u^m\,\L|\DD\mathsf P|^2}+\,2\,(1-m)\iring{u\,\mathsf P\,\DD\mathsf P\cdot\DD\L\mathsf P}\\
&&\hspace*{1cm}+\,2\,(1-m)\iring{u\,\DD\mathsf P\cdot\DD\mathsf P\,\L\mathsf P}-\,2\iring{u\,\DD\mathsf P\cdot\DD|\DD\mathsf P|^2}\\
&&+\,\alpha^2\,\idB S{\((u^m)'\,|\DD\mathsf P|^2 - u^m\,\tfrac\partial{\partial s}\,(|\DD\mathsf P|^2)\)}\\
&&\quad-\,\alpha^2\,\idB s{\((u^m)'\,|\DD\mathsf P|^2-u^m\,\tfrac\partial{\partial s}\,(|\DD\mathsf P|^2)\)}\\
&=&-\iring{u^m\,\L|\DD\mathsf P|^2}+\,2\,(1-m)\iring{u\,\mathsf P\,\DD\mathsf P\cdot\DD\L\mathsf P}\\
&&+\,2\,(1-m)\iring{u\,\DD\mathsf P\cdot\DD\mathsf P\,\L\mathsf P}\\
&&+\,\alpha^2\,\idB S{\((u^m)'\,|\DD\mathsf P|^2+u^m\,\tfrac\partial{\partial s}\,(|\DD\mathsf P|^2)\)}\\
&&\quad-\,\alpha^2\,\idB s{\((u^m)'\,|\DD\mathsf P|^2+u^m\,\tfrac\partial{\partial s}\,(|\DD\mathsf P|^2)\)}\,,
\end{eqnarray*}
where the last line is given by an integration by parts, upon exploiting the identity $u\,\DD\mathsf P=-\,\DD(u^m)$:
\begin{multline*}
\iring{u\,\DD\mathsf P\cdot\DD|\DD\mathsf P|^2}=-\iring{\DD(u^m)\cdot\DD|\DD\mathsf P|^2}\\
=\iring{u^m\,\L|\DD\mathsf P|^2}
-\,\alpha^2\,\idB S{u^m\,\tfrac\partial{\partial s}\,(|\DD\mathsf P|^2)}
+\alpha^2\,\idB s{u^m\,\tfrac\partial{\partial s}\,(|\DD\mathsf P|^2)}\,.
\end{multline*}
1) Using the definition of $\mathcal A[\mathsf P]$, we get that
\be{Id1}
-\iring{u^m\,\L|\DD\mathsf P|^2}=-\,2\iring{u^m\,\mathcal A[\mathsf P]}-\,2\iring{u^m\,\DD\mathsf P\cdot\DD\L\mathsf P}\,.
\ee
2) Taking advantage again of $u\,\DD\mathsf P=-\,\DD(u^m)$, an integration by parts gives
\begin{multline*}
\hspace*{-12pt}\iring{u\,\DD\mathsf P\cdot\DD\mathsf P\,\L\mathsf P}=-\iring{\DD(u^m)\cdot\DD\mathsf P\,\L\mathsf P}\\
=\iring{u^m\,(\L\mathsf P)^2}\,+\iring{u^m\,\DD\mathsf P\cdot\DD\L\mathsf P}\\
-\,\alpha^2\,\idB S{u^m\,\mathsf P'\L\mathsf P}+\alpha^2\,\idB s{u^m\,\mathsf P'\L\mathsf P}\,.
\end{multline*}
and, with $u\,\mathsf P=\frac m{1-m}\,u^m$, we find that
\begin{multline}
2\,(1-m)\iring{u\,\mathsf P\,\DD\mathsf P\cdot\DD\L\mathsf P}+\,2\,(1-m)\iring{u\,\DD\mathsf P\cdot\DD\mathsf P\,\L\mathsf P}\\
=\,2\,(1-m)\iring{u^m\,(\L\mathsf P)^2}+2\iring{u^m\,\DD\mathsf P\cdot\DD\L\mathsf P}\hspace*{4cm}\\
-\,2\,(1-m)\,\alpha^2\,\idB S{u^m\,\mathsf P'\L\mathsf P}+2\,(1-m)\,\alpha^2\,\idB s{u^m\,\mathsf P'\L\mathsf P}\,.\label{Id2}
\end{multline}
Summing~\eqref{Id1} and~\eqref{Id2}, using~\eqref{b} and passing to the limits as $s\to0_+$, $S\to+\infty$, establishes~\eqref{BLWL}.\end{proof}

\subsection{Second step: two remarkable identities.} Let us define
\[
\mathsf k[\mathsf P]:=\tfrac12\,\Delta_\omega\,|\nabla_\omega\mathsf P|^2-\nabla_\omega\mathsf P\cdot\nabla_\omega\Delta_\omega\,\mathsf P-\tfrac1{n-1}\,(\Delta_\omega\,\mathsf P)^2-(n-2)\,\alpha^2\,|\nabla_\omega\mathsf P|^2
\]
and
\[
\mathcal R[\mathsf P]:=\mathcal K[\mathsf P]-\(\frac1n-(1-m)\)\(\L\mathsf P\)^2\,.
\]
We observe that
\[
\mathcal R[\mathsf P]=\frac12\,\L\,|\DD\mathsf P|^2-\,\DD\mathsf P\cdot\DD\L\mathsf P-\frac1n\,(\L\mathsf P)^2
\]
is independent of $m$. We recall the result of~\cite[Lemma~5.1]{DEL2015} and give its proof for completeness.
\par\smallskip\begin{lemma}\label{Lem:Derivmatrixform1}{\sl Let $d\in\N$, $n\in\R$ such that $n>d\ge2$, and consider a function $\mathsf P\in C^3(\R^d\setminus\{0\})$. Then,
\[
\mathcal R[\mathsf P]=\alpha^4\(1-\frac1n\)\left[\mathsf P''-\frac{\mathsf P'}s-\frac{\Delta_\omega\,\mathsf P}{\alpha^2\,(n-1)\,s^2}\right]^2+\frac{2\,\alpha^2}{s^2}\left|\nabla_\omega\mathsf P'-\frac{\nabla_\omega\mathsf P}s \right|^2+\frac{\mathsf k[\mathsf P]}{s^4}\,.
\]}\end{lemma}\par\smallskip
\begin{proof} By definition of $\mathcal R[\mathsf P]$, we have
\begin{eqnarray*}
\mathcal R[\mathsf P]&=&\frac{\alpha^2}2\left[\alpha^2\,\mathsf P'^2+\frac{|\nabla_\omega\mathsf P|^2}{s^2}\right]''+\frac{\alpha^2}2\frac{n-1}s
\left[\alpha^2\,\mathsf P'^2+\frac{|\nabla_\omega\mathsf P|^2}{s^2}\right]'+\frac1{2\,s^2}\,\Delta_\omega\left[\alpha^2\,\mathsf P'^2+\frac{|\nabla_\omega\mathsf P|^2}{s^2}\right]\\
&&-\,\alpha^2\,\mathsf P'\(\alpha^2\,\mathsf P''+\alpha^2\,\frac{n-1}s\,\mathsf P'+\frac{\Delta_\omega\,\mathsf P}{s^2}\)'-\frac1{s^2}
\nabla_\omega\mathsf P\cdot\nabla_\omega\(\alpha^2\,\mathsf P''+\alpha^2\,\frac{n-1}s\,\mathsf P'+\frac{\Delta_\omega\,\mathsf P}{s^2}\)\\
&&-\,\frac 1n\(\alpha^2\,\mathsf P''+\alpha^2\,\frac{n-1}s\,\mathsf P'+\frac{\Delta_\omega\,\mathsf P}{s^2}\)^2\,,
\end{eqnarray*}
which can be expanded as
\begin{eqnarray*}
\mathcal R[\mathsf P]&=&\frac{\alpha^2}2\left[ 2\,\alpha^2\,\mathsf P''^2+2\,\alpha^2\,\mathsf P'\,\mathsf P'''+2\,\frac{|\nabla_\omega\mathsf P'|^2+\nabla_\omega\mathsf P\cdot\nabla_\omega\mathsf P''}{s^2}
-8\,\frac{\nabla_\omega\mathsf P\cdot\nabla_\omega\mathsf P'}{s^3}+6\,\frac{|\nabla_\omega\mathsf P|^2}{s^4}\right]\\
&&+\,\alpha^2\,\frac{n-1}s\left[\alpha^2\,\mathsf P'\,\mathsf P''+\frac{\nabla_\omega\mathsf P\cdot\nabla_\omega\mathsf P'}{s^2}-\frac{|\nabla_\omega\mathsf P|^2}{s^3}\right]+\frac1{s^2}\left[\alpha^2\,\mathsf P'\Delta_\omega\,\mathsf P'+\alpha^2\,|\nabla_\omega\mathsf P'|^2+\frac{\Delta_\omega\,|\nabla_\omega\mathsf P|^2}{2\,s^2}\right]\\
&&-\,\alpha^2\,\mathsf P'\(\alpha^2\,\mathsf P'''+\alpha^2\,\frac{n-1}s\,\mathsf P''-\,\alpha^2\,\frac{n-1}{s^2}\mathsf P'-2\,\frac{\Delta_\omega\,\mathsf P}{s^3}+\frac{\Delta_\omega\,\mathsf P'}{s^2}\)\\
&&\hspace*{2cm}-\frac1{s^2}
\(\alpha^2\,\nabla_\omega\mathsf P\cdot\nabla_\omega\mathsf P''+\alpha^2\,\frac{n-1}s\nabla_\omega\mathsf P\cdot\nabla_\omega\mathsf P'+\frac{\nabla_\omega\mathsf P\cdot\nabla_\omega\Delta_\omega\,\mathsf P}{s^2}\)\\
&&-\,\frac 1n\left[\alpha^4\,\mathsf P''^2+\alpha^4\,\frac{(n-1)^2}{s^2}\,\mathsf P'^2+\frac{(\Delta_\omega\,\mathsf P)^2}{s^4}+2\,\alpha^4\,\frac{n-1}s\,\mathsf P'\,\mathsf P''+2\,\alpha^2\,\frac{\mathsf P''\Delta_\omega\,\mathsf P}{s^2}+2\,\alpha^2\,\frac{n-1}{s^3}\mathsf P'\Delta_\omega\,\mathsf P\right]\,.
\end{eqnarray*}
Collecting terms proves the result.\end{proof}

\medskip Now let us study the quantity $\mathsf k[\mathsf P]$ which appears in the statement of Lemma~\ref{Lem:Derivmatrixform1}.
The following computations are adapted from~\cite{MR3229793} and~\cite[Section~5]{DEL2015}. For completeness, we give a simplified proof in the special case of the sphere $(\S^{d-1},g)$ considered as a Riemannian manifold with standard metric~$g$. We denote by $\mathrm Hf$ the \emph{Hessian} of~$f$, which is seen as $(d-1)\times(d-1)$ matrix, identify its trace with the Laplace-Beltrami operator on $\S^{d-1}$ and use the notation $\|\mathrm A\|^2:=\mathrm A:\mathrm A$ for the sum of the squares of the coefficients of the matrix~$A$. It is convenient to define the \emph{trace free Hessian}, the tensor $\mathrm Zf$ and its trace free counterpart respectively~by
\[
\mathrm Lf:=\mathrm Hf-\frac1{d-1}\,(\Delta_\omega f)\,g\,,\quad\mathrm Zf:=\frac{\nabla_\omega f\otimes\nabla_\omega f}f\quad\mbox{and}\quad\mathrm Mf:=\mathrm Zf-\frac1{d-1}\,\frac{|\nabla_\omega f|^2}f\,g
\]
whenever $f\neq0$. Elementary computations show that
\be{TraceFree}
\|\mathrm Lf\|^2=\|\mathrm Hf\|^2-\frac1{d-1}\,(\Delta_\omega f)^2\quad\mbox{and}\quad\|\mathrm Mf\|^2=\|\mathrm Zf\|^2-\frac1{d-1}\,\frac{|\nabla_\omega f|^4}{f^2}=\frac{d-2}{d-1}\,\frac{|\nabla_\omega f|^4}{f^2}\,.
\ee
The Bochner-Lichnerowicz-Weitzenb\"ock formula on $\S^{d-1}$ takes the simple form
\be{BLW}
\tfrac12\,\Delta_\omega\,(|\nabla_\omega f|^2)=\|\mathrm Hf\|^2+\nabla_\omega(\Delta_\omega f)\cdot\nabla_\omega f+(d-2)\,|\nabla_\omega f|^2
\ee
where the last term, \emph{i.e.}, $ \mathrm{Ric}(\nabla_\omega f,\nabla_\omega f)=(d-2)\,|\nabla_\omega f|^2$, accounts for the Ricci curvature tensor contracted with \hbox{$\nabla_\omega f \otimes\nabla_\omega f$}.

We recall that $\alpha_{\rm FS}:=\sqrt{\frac{d-1}{n-1}}$ and $\nu=1/(1-m)$. Let us introduce the notations
\[
\delta:=\frac1{d-1}-\frac1{n-1}
\]
and
\[
\mathcal B[\mathsf P]:=\isph{\(\tfrac12\,\Delta_\omega(|\nabla_\omega\mathsf P|^2)-\nabla_\omega(\Delta_\omega\mathsf P)\cdot\nabla_\omega\mathsf P-\tfrac1{n-1}\,(\Delta_\omega\mathsf P)^2\)\,\mathsf P^{1-\nu}}\,,
\]
so that
\[
\isph{\mathsf k[\mathsf P]\,\mathsf P^{1-\nu}}=\mathcal B[\mathsf P]-(n-2)\,\alpha^2\isph{|\nabla_\omega\mathsf P|^2\,\mathsf P^{1-\nu}}\,.
\]
\par\smallskip\begin{lemma}\label{kappapositive}{\sl Assume that $d\ge2$ and $1/(1-m)=\nu>n>d$. There exists a positive constant $c(n,m,d)$ such that, for any positive function $\mathsf P\in C^3({\S^{d-1}})$,
\[
\isph{\mathsf k[\mathsf P]\,\mathsf P^{1-\nu}}\ge(n-2)\,\big(\alpha_{\rm FS}^2-\,\alpha^2\big)\isph{|\nabla_\omega\mathsf P|^2\,\mathsf P^{1-\nu}}+c(n,m,d)\isph{\frac{|\nabla_\omega\mathsf P|^4}{\mathsf P^2}\,\mathsf P^{1-\nu}}\,.
\]}\end{lemma}\par\smallskip
\begin{proof} If $d=2$, we identify $\S^1$ with $[0,2\pi)\ni\theta$ and denote by $\mathsf P_\theta$ and $\mathsf P_{\theta\theta}$ the first and second derivatives of $\mathsf P$ with respect to $\theta$. As in~\cite[Lemma~5.3]{DEL2015}, a direct computation shows that
\[
\mathsf k[\mathsf P]=\frac{n-2}{n-1}\,|\mathsf P_{\theta\theta}|^2-(n-2)\,\alpha^2\,|\mathsf P_\theta|^2=(n-2)\,\(\alpha_{\rm FS}^2\,|\mathsf P_{\theta\theta}|^2-\,\alpha^2\,\,|\mathsf P_\theta|^2\)\,.
\]
By the Poincar\'e inequality, we have
\[
\icircle{\left|\frac\partial{\partial\theta}\(\mathsf P^\frac{1-\nu}2\,\mathsf P_\theta\)\right|^2}\ge\icircle{\left|\mathsf P^\frac{1-\nu}2\,\mathsf P_\theta\right|^2}\,.
\]
On the other hand, an integration by parts shows that
\[
\icircle{\mathsf P_{\theta\theta}\,\frac{|\mathsf P_\theta|^2}{\mathsf P}\,\mathsf P^{1-\nu}}=\frac 13\icircle{\frac\partial{\partial\theta}\(|\mathsf P_\theta|^2\,\mathsf P_\theta\)\,\mathsf P^{-\nu}}=\frac\nu3\icircle{\frac{|\mathsf P_\theta|^4}{\mathsf P^2}\,\mathsf P^{1-\nu}}
\]
and, as a consequence, by expanding the square, we obtain
\[
\icircle{\left|\frac\partial{\partial\theta}\(\mathsf P^\frac{1-\nu}2\,\mathsf P_\theta\)\right|^2}=\icircle{\left|\mathsf P_{\theta\theta}+\frac{1-\nu}2\,\frac{|\mathsf P_\theta|^2}{\mathsf P}\right|^2\,\mathsf P^{1-\nu}}=\icircle{|\mathsf P_{\theta\theta}|^2\,\mathsf P^{1-\nu}}-\,\frac{(\nu-1)\,(\nu+3)}{12}\icircle{\frac{|\mathsf P_\theta|^4}{\mathsf P^2}\,\mathsf P^{1-\nu}}\,.
\]
The result follows with $c(n,m,2)=\frac{n-2}{n-1}\,\frac1{12}\,(\nu-1)\,(\nu+3)=\frac{n-2}{n-1}\,\frac{m\,(4-3\,m)}{12\,(1-m)^2}$ from
\[\label{Pincd2}
\icircle{|\mathsf P_{\theta\theta}|^2\,\mathsf P^{1-\nu}}\ge\icircle{|\mathsf P_\theta|^2\,\mathsf P^{1-\nu}}+\frac{(\nu-1)\,(\nu+3)}{12}\icircle{\frac{|\mathsf P_\theta|^4}{\mathsf P^2}\,\mathsf P^{1-\nu}}\,.
\]

Assume next that $d\ge3$. We follow the method of~\cite[Lemma~5.2]{DEL2015}. Applying~\eqref{BLW} with $f=\mathsf P$ and multiplying by $\mathsf P^{1-\nu}$ yields, after an integration on $\S^{d-1}$, that $\mathcal B[\mathsf P]$ can also be written as
\[
\mathcal B[\mathsf P]=\isph{ \(\|\mathrm H\mathsf P\|^2+(d-2)\,|\nabla_\omega\mathsf P|^2-\tfrac1{n-1}\,(\Delta_\omega\mathsf P)^2\)\,\mathsf P^{1-\nu}}\,.
\]
We recall that $n>d\ge3$ and set $\mathsf P=f^\beta$ with $\beta=\frac2{3-\nu}$. A straightforward computation shows that $\mathrm H f^\beta=\beta\,f^{\beta-1}\,\big(\mathrm H f+(\beta-1)\,\mathrm Z f\big)$ and hence
\begin{multline*}
\mathcal B[\mathsf P]=\beta^2 \isph{ \(\|\mathrm H f+(\beta-1)\,\mathrm Z f\|^2+(d-2)\,|\nabla_\omega f|^2-\tfrac1{n-1}\,\big(\mathrm{Tr}\,(\mathrm H f+(\beta-1)\,\mathrm Z f)\big)^2\)}\\
=\beta^2 \isph{ \(\|\mathrm L f+(\beta-1)\,\mathrm M f\|^2+(d-2)\,|\nabla_\omega f|^2+\delta\,\big(\mathrm{Tr}\,(\mathrm H f+(\beta-1)\,\mathrm Z f)\big)^2\)}\,.
\end{multline*}
Using~\eqref{TraceFree}, we deduce from
\begin{multline*}
\isph{\Delta_\omega\,f\,\frac{|\nabla_\omega f|^2}{f}}=\isph{ \frac{|\nabla_\omega f|^4}{f^2} }-2 \isph{\mathrm H f:\mathrm Z f}\\
=\frac{d-1}{d-2}\,\isph{\|\mathrm Mf\|^2}-2 \isph{\mathrm L f:\mathrm Z f}-\frac{2}{d-1} \isph{ \Delta_\omega\,f\,\frac{|\nabla_\omega f|^2}f}
\end{multline*}
that
\begin{multline*}
\isph{\Delta_\omega\,f\,\frac{|\nabla_\omega f|^2}{f}}=\frac{d-1}{d+1}\left[\isph{\frac{d-1}{d-2}\,\|\mathrm Mf\|^2}-2 \isph{\mathrm L f:\mathrm Z f} \right]\\
=\frac{d-1}{d+1}\left[\isph{\frac{d-1}{d-2}\,\|\mathrm Mf\|^2}-2 \isph{\mathrm L f:\mathrm M f} \right]
\end{multline*}
on the one hand, and from~\eqref{BLW} integrated on $\S^{d-1}$ that
\[
\isph{(\Delta_\omega\,f)^2}=\frac{d-1}{d-2} \isph{\|\mathrm L f\|^2}+(d-1)\isph{|\nabla_\omega f|^2}
\]
on the other hand. Hence we find that
\begin{multline*}
\isph{\big(\mathrm{Tr}\,(\mathrm H f+(\beta-1)\,\mathrm Z f)\big)^2}=\isph{\((\Delta_\omega\,f)^2+2\,(\beta-1)\,\Delta_\omega\,f\,\frac{|\nabla_\omega f|^2}{f}+(\beta-1)^2\,\frac{|\nabla_\omega f|^4}{f^2}\)}\\
\hspace*{1cm}=\frac{d-1}{d-2} \isph{\|\mathrm L f\|^2}+(d-1)\isph{|\nabla_\omega f|^2}\\
\hspace*{6cm}+2\,(\beta-1)\,\frac{d-1}{d+1}\left[\isph{\frac{d-1}{d-2}\,\|\mathrm Mf\|^2}-2 \isph{\mathrm L f:\mathrm M f} \right]\\
+(\beta-1)^2\,\frac{d-1}{d-2}\isph{\|\mathrm Mf\|^2}\,.
\end{multline*}
Altogether, we obtain
\[
\mathcal B[\mathsf P]=\beta^2\isph{\Big(\mathsf a\,\|\mathrm L f\|^2+\,2\,\mathsf b\,\mathrm L f:\mathrm M f+\,\mathsf c\,\|\mathrm M f\|^2\Big)}+\beta^2\,\big(d-2+\delta\,(d-1)\big)\isph{|\nabla_\omega f|^2}
\]
where
\[
\mathsf a=1+\delta\,\frac{d-1}{d-2}\,,\quad\mathsf b=(\beta-1)\,\(1-\,2\,\delta\,\frac{d-1}{d+1}\)\quad\mbox{and}\quad\mathsf c=(\beta-1)^2\,\(1+\delta\,\frac{d-1}{d-2}\)+2\,(\beta-1)\,\frac{\delta\,(d-1)^2}{(d+1)\,(d-2)}\,.
\]
A tedious but elementary computation shows that
\[
\mathcal B[\mathsf P]=\mathsf a\,\beta^2\isph{\left\|\mathrm L f+\tfrac{\mathsf b}{\mathsf a}\;\mathrm M f\right\|^2}+\big(\mathsf c-\tfrac{\mathsf b^2}{\mathsf a}\big)\,\beta^2\isph{\left\|\mathrm M\,f\right\|^2}+\beta^2\,(n-2)\,\alpha_{\rm FS}^2\isph{|\nabla_\omega f|^2}
\]
can be written in terms of $ \mathsf P$ as
\[
\mathcal B[\mathsf P]=\isph{\mathrm Q[\mathsf P]\,\mathsf P^{1-\nu}}+(n-2)\,\alpha_{\rm FS}^2\isph{|\nabla_\omega\mathsf P|^2\,\mathsf P^{1-\nu}}
\]
where
\[
\mathrm Q[\mathsf P]:=\alpha_{\rm FS}^2\,\frac{n-2}{d-2}\,{\left\|\mathrm L\mathsf P+\tfrac{3\,(\nu-1)\,(n-d)}{(d+1)\,(n-2)\,(\nu-3)}\;\mathrm M\mathsf P\right\|^2}+\tfrac{(d-1)\,(\nu-1)\,(n-d)\,[((4\,d-5)\,n+d-8)\,\nu+9\,(n-d))]}{(d-2)\,(d+1)^2\,(\nu-3)^2\,(n-2)\,(n-1)}\,\|\mathrm M\mathsf P\|^2
\]
is positive definite. This concludes the proof in the case $d\ge3$ with $c(n,m,d)=\frac{m\,(n-d)\,[4\,(d+1)\,(n-2)-9\,m\,(n-d)]}{(d+1)^2\,(3\,m-2)^2\,(n-2)\,(n-1)}$.\end{proof}

Let us recall that
\[
\mathcal K[\mathsf P]=\mathcal R[\mathsf P]+\(\frac1n-(1-m)\)\(\L\mathsf P\)^2\,.
\]
We can collect the two results of Lemmas~\ref{Lem:Derivmatrixform1} and~\ref{kappapositive} as follows.
\par\smallskip\begin{corollary}\label{Cor:TwoIdentities}{\sl Let $d\in\N$, $n\in\R$ be such that $n>d\ge2$, and consider a positive function $\mathsf P\in C^3(\R^d\setminus\{0\})$. If $u$ is related to $\mathsf P$ by $\mathsf P=\frac m{1-m}\,u^{m-1}$ for some $m\in(1-\frac1n,1)$, then there exists a positive constant $c(n,m,d)$ such that 
\begin{multline*}
\irdmu{\mathcal R[\mathsf P]\,u^m}\ge\alpha^4\(1-\frac1n\)\irdmu{\left[\mathsf P''-\frac{\mathsf P'}s-\frac{\Delta_\omega\,\mathsf P}{\alpha^2\,(n-1)\,s^2}\right]^2\,u^m}+2\,\alpha^2\irdmu{\frac1{s^2}\,\left|\nabla_\omega\mathsf P'-\frac{\nabla_\omega\mathsf P}s \right|^2\,u^m}\\
+(n-2)\,\big(\alpha_{\rm FS}^2-\,\alpha^2\big)\irdmu{\frac1{s^4}\,|\nabla_\omega\mathsf P|^2\,u^m} + c(n,m,d)\irdmu{\frac1{s^4}\,\frac{|\nabla_\omega\mathsf P|^4}{\mathsf P^2}\,u^m} \,.
\end{multline*}}\end{corollary}\par\smallskip

\subsection{Third step: concavity of the R\'enyi entropy powers and consequences}

We keep investigating the properties of the flow defined by~\eqref{ordinaryflow}. Let us define the \emph{entropy} as
\[
\mathcal E:=\irdmu{u^m}
\]
and observe that
\[
\mathcal E'=(1-m)\,\mathcal I
\]
if $u$ solves~\eqref{FDE}, after integrating by parts. The fact that boundary terms do not contribute, \emph{i.e.},
\be{FisherBC}
\lim_{s\to0_+}\idB s{u^m\,\mathsf P'}=\lim_{S\to+\infty}\idB S{u^m\,\mathsf P'}=0
\ee
will be justified in Section~\ref{Sec:RegDecay}: see Proposition~\ref{Prop:b}. Note that we use $ ' $ both for derivation w.r.t.~$t$ and w.r.t.~$s$, at least when this does not create any ambiguity. As in Section~\ref{RenyiStrategy}, we introduce the \emph{R\'enyi entropy power}
\[
\mathcal F:=\mathcal E^\sigma
\]
for some exponent $\sigma$ to be chosen later, and find that $\mathcal F'=\sigma\,(1-m)\,\mathcal G$ where $\mathcal G:=\mathcal E^{\sigma-1}\,\mathcal I$. With $\mathcal H:=\mathcal E^{-\sigma}\,\mathcal G'$, by using Lemma~\ref{Lem:DerivFisherL}, we also find that $\mathcal E^{-\sigma}\,\mathcal F''=\sigma\,(1-m)\,\mathcal H$ where
\begin{multline*}
\mathcal E^2\,\mathcal H=\mathcal E^{2-\sigma}\,\mathcal G'=\frac1{\sigma\,(1-m)}\,\mathcal E^{2-\sigma}\,\mathcal F''=(1-m)\,(\sigma-1)\(\irdmu{u\,|\DD\mathsf P|^2}\)^2-\,2\irdmu{u^m}\irdmu{\mathcal K[\mathsf P]\,u^m}\\
=(1-m)\,(\sigma-1)\(\irdmu{u\,|\DD\mathsf P|^2}\)^2-\,2\,\(\frac1n-(1-m)\)\irdmu{u^m}\irdmu{\(\L\mathsf P\)^2\,u^m}\hspace*{-2cm}\\
-\,2\irdmu{u^m}\irdmu{\mathcal R[\mathsf P]\,u^m}
\end{multline*}
if $\lim_{s\to0_+}\mathsf b(s)=\lim_{S\to+\infty}\mathsf b(S)=0$. Using $u\,\DD\mathsf P=-\,\DD(u^m)$, we know that
\[
\irdmu{u\,|\DD\mathsf P|^2}=-\irdmu{\DD(u^m)\cdot\DD\mathsf P}=\irdmu{u^m\,\L\mathsf P}
\]
and so, with the choice
\[
\sigma=\frac2n\,\frac1{1-m}-1\,,
\]
we may argue as in Section~\ref{RenyiStrategy} and get that
\[
\mathcal E^2\,\mathcal H+(1-m)\,(\sigma-1)\,\mathcal E \irdmu{u^m\,\left|\L\mathsf P - \frac{\irdmu{u\,|\DD\mathsf P|^2}}{\irdmu{u^m}}\right|^2}+\,2\,\mathcal E \irdmu{\mathcal R[\mathsf P]\,u^m}=0
\]
if $\lim_{s\to0_+}\mathsf b(s)=\lim_{S\to+\infty}\mathsf b(S)=0$. So, if $\alpha\le\alpha_{\rm FS}$ and $\mathsf P$ is of class $C^3$, by Corollary~\ref{Cor:TwoIdentities}, as a function of $t$, $\mathcal F$ is concave, that is, $\mathcal G=\mathcal E^{\sigma-1}\,\mathcal I$ is non-increasing in $t$. Formally, $\mathcal G$ converges towards a minimum, for which necessarily $\L\mathsf P$ is a constant and $\mathcal R[\mathsf P]=0$, which proves that $\mathsf P(x)=\mathsf a+\mathsf b\,|x|^2$ for some real constants $\mathsf a$ and~$\mathsf b$, according to Corollary~\ref{Cor:TwoIdentities}. Since $\frac{2\,(1-\vartheta)}{\vartheta\,(p+1)}=\sigma-1$, the minimization of $\mathcal G$ under the mass constraint $\irdmu u=\irdmu{v^{2p}}$ is equivalent to the \emph{Caffarelli-Kohn-Nirenberg interpolation inequalities}~\eqref{CKN}, since for some constant~$\kappa$,
\[
\mathcal G=\mathcal E^{\sigma-1}\,\mathcal I=\kappa\,\(\irdmu{v^{p+1}}\)^{\sigma-1}\,\irdmu{|\DD v|^2}\quad\mbox{with}\quad v=u^{m-1/2}\,.
\]
We emphasize that~\eqref{FDE} preserves mass, that is, $\frac d{\dt}\irdmu{v^{2p}}=\frac d{\dt}\irdmu u=\irdmu{\L u^m}=0$ because, as we shall see in Proposition~\ref{Prop:b}, no boundary terms appear when integrating by parts if $v$ is an extremal function associated with~\eqref{CKN1}. In particular, for mass conservation we need
\be{MassBC}
\lim_{s\to0_+}\idB s{u\,\mathsf P'}=\lim_{S\to+\infty}\idB S{u\,\mathsf P'}=0 \, .
\ee

The above remarks on the monotonicity of $\mathcal G$ and the symmetry properties of its minimizers can in fact be extended to the analysis of the symmetry properties of all critical points of $\mathcal G$. This is actually the contents of Theorem~\ref{Thm:Rigidity}.

\medskip\noindent\emph{Proof of Theorem~\ref{Thm:Rigidity}.} Let $w$ be a positive solution of equation~\eqref{ELeq}. As pointed out above, by choosing
\[
w(x)=u^{m-1/2}(r^\alpha,\omega)\,,
\]
we know that $u$ is a critical point of $\mathcal G$ under a mass constraint on $\ird u$, so that we can write the corresponding Euler-Lagrange equation as $\mathrm d\mathcal G[u]=C$, for some constant $C$. That is, $\irdmu{\mathrm d\mathcal G[u]\cdot\L u^m }= C\,\irdmu{\L u^m }=0$ thanks to~\eqref{MassBC}. Using $\L u^m$ as a test function amounts to apply the flow of~\eqref{FDE} to $\mathcal G$ with initial datum $u$ and compute the derivative with respect to $t$ at $t=0$. This means
\begin{multline*}
0=\irdmu{\mathrm d\mathcal G[u]\cdot\L u^m}=\mathcal E^\sigma\,\mathcal H\\
=-\,(1-m)\,(\sigma-1)\,\mathcal E^{\sigma-1}\irdmu{u^m\,\left|\L\mathsf P-\frac{\irdmu{u\,|\DD\mathsf P|^2}}{\irdmu{u^m}}\right|^2}-\,2\,\mathcal E^{\sigma-1}\irdmu{\mathcal R[\mathsf P]\,u^m}
\end{multline*}
if $\lim_{s\to0_+}\mathsf b(s)=\lim_{S\to+\infty}\mathsf b(S)=0$ and~\eqref{FisherBC} holds. Here we have used Lemma~\ref{Lem:DerivFisherL}. We emphasize that this proof is purely variational and does not rely on the properties of the solutions to~\eqref{FDE}, although using the flow was very useful to explain our strategy. All we need is that no boundary term appears in the integrations by parts. Hence, in order to obtain a complete proof, we have to prove that~\eqref{FisherDerBC},~\eqref{FisherBC} and~\eqref{MassBC} hold with $\mathsf b$ defined by~\eqref{b}, whenever $u$ is a critical point of $\mathcal G$ under mass constraint. This will be done in Proposition~\ref{Prop:b}. Using Corollary~\ref{Cor:TwoIdentities}, we know that $\mathcal R[\mathsf P]=0$, $\nabla_\omega\mathsf P=0$ a.e.~in~$\R^d$ and $\L\mathsf P=\frac{\irdmu{u\,|\DD\mathsf P|^2}}{\irdmu{u^m}}$ a.e.~in~$\R^d$, with $\mathsf P=\frac m{1-m}\,u^{m-1}$. We conclude as in~\cite[Corollary 5.5]{DEL2015} that $\mathsf P$ is an affine function of~$s^2$.
\hfill\ \qed

\section{Regularity and decay estimates}\label{Sec:RegDecay}

In this last section we prove the regularity and decay estimates on $w$ (or on $\mathsf P$ or $u$) that are necessary to establish the absence of boundary terms in the integrations by parts of Section~\ref{Sec:BE}.
\par\smallskip\begin{proposition}\label{Prop:b}{\sl Under Condition~\eqref{parameters}, if $w$ is a positive solution in $\mathrm H^p_{\beta,\gamma}(\R^d)$ of~\eqref{ELeq}, then~\eqref{FisherDerBC},~\eqref{FisherBC} and~\eqref{MassBC} hold with $\mathsf b$ as defined by~\eqref{b}, $u=v^{2p}$ and $v$ given by~\eqref{wv}.}\end{proposition}\par\smallskip
To prove this result, we split the proof in several steps: we will first establish a uniform bound and a decay rate for $ w $ inspired by~\cite{DMN2015} in Lemmas~\ref{Lem:estimates1}, \ref{Lem:estimates2}, and then follow the methodology of~\cite[Appendix]{DEL2015} in the subsequent Lemma~\ref{Lem:decayinRd}.
\par\smallskip\begin{lemma}\label{Lem:estimates1}{\sl Let $\beta$, $\gamma$ and $p$ satisfy the relations~\eqref{parameters}. Any positive solution $w$ of~\eqref{ELeq} such that
\be{energyintegrals}
\nrm w{2p,\gamma}+\nrm{\nabla w}{2,\beta}+\nrm w{p+1,\gamma}^{1-\vartheta}<+\infty\,.
\ee
is uniformly bounded and tends to $0$ at infinity, uniformly in $|x|$.}\end{lemma}\par\smallskip
\begin{proof} The strategy of the first part of the proof is similar to the one in~\cite[Lemma 3.1]{DMN2015}, which was restricted to the case $\beta=0$.

Let us set $\delta_0:=2\,(p_\star-p)$. For any $A>0$, we multiply~\eqref{ELeq} by $(w\wedge A)^{1+\delta_0}$ and integrate by parts (or, equivalently, plug it in the weak formulation of~\eqref{ELeq}): we point out that the latter is indeed an admissible test function since $ w \in \mathrm H^p_{\beta,\gamma}(\R^d) $. In that way, by letting $ A \to+\infty $, we obtain the identity
\[\label{eq: prima-stima-H1}
\frac{4\,(1+\delta_0)}{(2+\delta_0)^2}\ird{\left|\nabla{w^{1+\delta_0/2}}\right|^2\,|x|^{-\beta}}+\ird{w^{p+1+\delta_0}\,|x|^{-\gamma}}=\ird{w^{2p+\delta_0}\,|x|^{-\gamma}}\,.
\]
By applying~\eqref{CKN} with $p=p_\star$ (so that $\vartheta=1$) to the function $w=w^{1+{\delta_0}/2}$, we deduce that
\[
\nrm{w}{2p+\delta_1,\gamma}^{2+\delta_0}\le\frac{(2+\delta_0)^2}{4\,(1+\delta_0)}\,\C_{\beta,\gamma,p_\star}^2\,\nrm{w}{2p+\delta_0,\gamma}^{2p+\delta_0}
\]
with $2\,p+\delta_1=p_\star\,(2+\delta_0)$. Let us define the sequence $\{\delta_n\}$ by the induction relation $\delta_{n+1}:=p_\star\,(2+\delta_n)-2\,p$ for any $n\in\N$, that is,
\[\label{eq: stima-Lq-weighted-rec-solved}
\textstyle\delta_n=2\,\frac{p_\star-p}{p_\star-1}\(p_\star^{n+1}-1\)\quad\forall\,n\in\N\,,
\]
and take $q_n=2\,p+\delta_n$. If we repeat the above estimates with $\delta_0$ replaced by $\delta_n$ and $\delta_1$ replaced by $\delta_{n+1}$, we get
\[
\nrm{w}{q_{n+1},\gamma}^{2+\delta_n}\le\frac{(2+\delta_n)^2}{4\,(1+\delta_n)}\,\C_{\beta,\gamma,p_\star}^2\,\nrm{w}{q_n,\gamma}^{q_n}\,.
\]
By iterating this estimate, we obtain the estimate
\[
\nrm{w}{q_n,\gamma}\le C_n\,\nrm{w}{2p_\star,\gamma}^{\zeta_n}\quad\mbox{with}\quad\zeta_n=\frac{(p_\star-1)\,p_\star^n}{p-1+(p_\star-p)\,p_\star^n}
\]
where the sequence $\{C_n\}$ is defined by $C_0=1$ and
\[
C_{n+1}^{2+\delta_n}=\frac{(2+\delta_n)^2}{4\,(1+\delta_n)}\,\C_{\beta,\gamma,p_\star}^2\,C_n^{q_n}\quad\forall\,n\in\N\,.
\]
The sequence $\{C_n\}$ converges to a finite limit $C_\infty$. Letting $n\to\infty$ we obtain the uniform bound
\[\label{eq: stima-infty-indip}
\nrm{w}\infty\le C_\infty\,\nrm{w}{2p_\star,\gamma}^{\zeta_\infty}\le C_\infty\(\C_{\beta,\gamma,p_\star}\,\nrm{\nabla w}{2,\beta}\)^{\zeta_\infty}\le C_\infty\(\C_{\beta,\gamma,p_\star}\,\nrm w{2p,\gamma}^p\)^{\zeta_\infty}
\]
where $\zeta_\infty=\frac{p_\star-1}{p_\star-p}=\lim_{n\to\infty}\zeta_n$.

\medskip In order to prove that $\lim_{|x|\to+\infty}w(x)=0$, we can suitably adapt the above strategy. We shall do it as follows: we truncate the solution so that the truncated function is supported outside of a ball of radius $R_0$ and apply the iteration scheme. Up to an enlargement of the ball, that is, outside of a ball of radius $R_\infty=a\,R_0$ for some fixed numerical constant $a>1$, we get that $\left\| w \right\|_{\mathrm L^{\infty}(B_{R_\infty}^c)}$ is bounded by the energy localized in $B_{R_0}^c$. The conclusion will hold by letting $ R_0\to+\infty$. Let us give some details.

Let $ \xi\in C^\infty(\R^+) $ be a cut-off function such that $ 0 \le \xi \le 1 $, $ \xi \equiv 0 $ in $ [0,1) $ and $ \xi \equiv 1 $ in $ (2,+\infty) $. Given $ R_0 \ge 1 $, consider the sequence of radii defined by
\[
R_{n+1}=\( 1+\frac1{n^2} \) R_n \quad \forall\,n\in \mathbb{N}\,.
\]
By taking logarithms, it is immediate to deduce that $ \{ R_n \} $ is monotone increasing and that there exists $ a>1 $ such that
\[\label{lim-Rn}
R_\infty:=\lim_{n \to \infty} R_n = a\,R_0\,.
\]
Let us then define the sequence of radial cut-off functions $ \{ \xi_n \} $ by
\[
\xi_n(x) := \xi^2\!\( \frac{|x|-R_n}{R_{n+1}-R_n}+1 \) \quad \forall\,x \in\R^d\,.
\]
Direct computations show that there exists some constant $c>0$, which is independent of $n$ and $R_0$, such that
\be{estimates-cutoff}
\left| \nabla \xi_n(x) \right| \le c\,\frac{n^2}{R_n}\,\chi_{B_{R_{n+1}} \setminus B_{R_{n}}}\,, \quad \left| \nabla \xi_n^{1/2}(x) \right| \le c\,\frac{n^2}{R_n}\,\chi_{B_{R_{n+1}} \setminus B_{R_{n}}}\,, \quad \left| \Delta \xi_n(x) \right| \le c\,\frac{n^4}{R_n^2}\,\chi_{B_{R_{n+1}} \setminus B_{R_{n}}} \quad \forall\,x \in\R^d\,.
\ee
{}From here on we denote by $c$, $c'$, \emph{etc.} positive constants which are all independent of $n$ and $R_0$. We now introduce the analogue of the sequence $ \{ \delta_n \} $ above, which we relabel $ \{ \sigma_n \} $ to avoid confusion. Namely, we set $ \sigma_0:=2\,p-2 $ and $ \sigma_{n+1}=p_\star\,(2+\sigma_n)-2 $, so that $ \sigma_n=2\,(p\,p_\star^{n}-1) $. If we multiply~\eqref{ELeq} by $ \xi_n w^{1+\sigma_n} $ and integrate by parts, we obtain:
\[
\ird{\nabla{\( \xi_n\,w^{1+\sigma_n} \)} \cdot \nabla{w}\,|x|^{-\beta}}+\ird{\xi_n\,w^{p+1+\sigma_n}\,|x|^{-\gamma}} =\ird{\xi_n\,w^{2p+\sigma_n}\,|x|^{-\gamma}}\,,
\]
whence
\[\label{estimates-cutoff-1}
\frac{4\,(1+\sigma_n)}{(2+\sigma_n)^2}\ird{\xi_n\,\left| \nabla w^{1+\sigma_n/2} \right|^2\,|x|^{-\beta}}+\frac1{2+\sigma_n}\ird{\nabla{\xi_n} \cdot \nabla{w^{2+\sigma_n}}\,|x|^{-\beta}} \le\int_{B_{R_n}^c} w^{2p+\sigma_n}\,|x|^{-\gamma}\,dx\,.
\]
By integrating by parts the second term in the l.h.s. and combining this estimate with
\[
\ird{\left| \nabla\( \xi_n^{1/2}\,w^{1+\sigma_n/2} \) \right|^2\,|x|^{-\beta}} \le 2\ird{\xi_n\,\left| \nabla w^{1+\sigma_n/2} \right|^2\,|x|^{-\beta}}+2\ird{\left| \nabla \xi_n^{1/2}\right|^2 w^{2+\sigma_n}\,|x|^{-\beta}}\,,
\]
we end up with
\[\label{estimates-cutoff-2}
\begin{aligned}
\frac{2\,(1+\sigma_n)}{(2+\sigma_n)^2}\,\ird{\left| \nabla\( \xi_n^{1/2}\,w^{1+\sigma_n/2} \) \right|^2\,|x|^{-\beta}} - \frac{4\,(1+\sigma_n)}{(2+\sigma_n)^2}\,\ird{\left| \nabla \xi_n^{1/2} \right|^2 w^{2+\sigma_n}\,|x|^{-\beta}} & \\
-\,\frac1{2+\sigma_n}\,\ird{\( |x|^{-\beta} \Delta \xi_n - \beta\,|x|^{-\beta-2} x \cdot \nabla{\xi_n} \) w^{2+\sigma_n}} & \le\int_{B_{R_n}^c} w^{2p+\sigma_n}\,|x|^{-\gamma}\,dx\,.
\end{aligned}
\]
Thanks to~\eqref{estimates-cutoff}, we can deduce that
\begin{multline*}
\ird{\left|\nabla\(\xi_n^{1/2}\,w^{1+\sigma_n/2}\)\right|^2\,|x|^{-\beta}}\le\int_{B_{R_{n+1}}\setminus B_{R_{n}}}\(\frac{2\,c^2+c}{R_n^2}\,n^4+\frac{\beta\,c}{R_n}\,n^2\,|x|^{-1}\)w^{2+\sigma_n}\,|x|^{-\beta}\,dx\\
+\frac{(2+\sigma_n)^2}{2\,(1+\sigma_n)}\,\int_{B_{R_n}^c} w^{2p+\sigma_n}\,|x|^{-\gamma}\,dx\,.
\end{multline*}
In particular,
\[\label{estimates-cutoff-4}
\ird{\left| \nabla\( \xi_n^{1/2}\,w^{1+\sigma_n/2} \) \right|^2\,|x|^{-\beta}} \le c^\prime n^4\,\int_{B_{R_n}^c} w^{2+\sigma_n}\,|x|^{-\beta-2}\,dx+\frac{(2+\sigma_n)^2}{2\,(1+\sigma_n)}\,\| w \|_\infty^{2p-2}\,\int_{B_{R_n}^c} w^{2+\sigma_n}\,|x|^{-\gamma}\,dx\,.
\]
Since~\eqref{parameters} implies that $ \beta+2>\gamma $, by exploiting the explicit expression of $\sigma_n$ and applying~\eqref{CKN} with $p=p_\star$ (and $\vartheta=1$) to the function $ \xi_n^{1/2}\,w^{1+\sigma_n/2} $, we can rewrite our estimate as
\[\label{estimates-cutoff-5}
\left\| w \right\|^{2+\sigma_n}_{\mathrm L^{2+\sigma_{n+1},\gamma}(B_{R_{n+1}}^c)} \le c^{\prime\prime} p_\star^n \left\| w \right\|^{2+\sigma_n}_{\mathrm L^{2+\sigma_n,\gamma}(B_{R_n}^c)}\,.
\]
After iterating the scheme and letting $ n \to \infty $ we end up with
\[
\left\| w \right\|_{\mathrm L^{\infty}(B_{R_\infty }^c)} \le c^{\prime\prime\prime} \left\| w \right\|_{\mathrm L^{2p,\gamma}(B_{R_0}^c)}\,.
\]
Since $w$ is bounded in $\mathrm L^{2p,\gamma}(\R^d) $, in order to prove the claim it is enough to let $ R_0 \to+\infty $.
\end{proof}

\par\smallskip\begin{lemma}\label{Lem:estimates2}{\sl Let $\beta$, $\gamma$ and $p$ satisfy the relations~\eqref{parameters}. Any positive solution $w$ of~\eqref{ELeq} satisfying
\eqref{energyintegrals} is such that $w\in C^\infty(\R^d\setminus\{0\})$ and there exists two positive constants, $C_1$ and $C_2$ with $C_1<C_2$, such that for $|x|$ large enough,
\[
C_1\,|x|^{(\gamma-2-\beta)/(p-1)}\le w(x)\le C_2\,|x|^{(\gamma-2-\beta)/(p-1)}\,.
\]}\end{lemma}\par\smallskip
\begin{proof} By Lemma~\ref{Lem:estimates1} and elliptic bootstrapping methods we know that $w\in C^\infty(\R^d\setminus\{0\})$. Let us now consider the function $h(x):=C\,|x|^{(\gamma-2-\beta)/(p-1)}$, which satisfies the differential inequality
\[
-\,\mbox{div}\,\(|x|^{-\beta}\,\nabla h\)+(1-\varepsilon)\,|x|^{-\gamma}\,h^p\ge0\quad\forall\,x\in\R^d\setminus\{0\}
\]
for any $\varepsilon\in(0,1)$ and $C$ such that $C^{p-1}>\frac{2-\gamma+\beta}{1-\varepsilon}\,\frac{d-\gamma-p\,(d-2-\beta)}{(p-1)^2}$. On the other hand, by Lemma~\ref{Lem:estimates1}, $w^{2p-1}$ is negligible compared to $w^p$ as $|x|\to\infty$ and, as a consequence, for any $\varepsilon>0$ small enough, there is an $R_\varepsilon>0$ such that
\[
-\,\mbox{div}\,\(|x|^{-\beta}\,\nabla w\)+(1-\varepsilon)\,|x|^{-\gamma}\,w^p\le0\quad\mbox{if}\quad|x|\ge R_\varepsilon\,.
\]
With $q:=(1-\varepsilon)\,|x|^{-\gamma}\,\frac{h^p-w^p}{h-w}\ge0$, it follows that
\[
-\,\mbox{div}\,\(|x|^{-\beta}\,\nabla (h-w)\)+q\,(h-w)\ge0\quad\mbox{if}\quad|x|\ge R_\varepsilon\,.
\]
Hence, for $C$ large enough, we know that $h(x)\ge w(x)$ for any $x\in\R^d$ such that $|x|=R_\varepsilon$, and we also have that $\lim_{|x|\to+\infty}\big(h(x)-w(x)\big)=0$. Using the Maximum Principle, we conclude that $0\le w(x)\le h(x)$ for any $x\in\R^d$ such that $|x|\ge R_\varepsilon$. The lower bound uses a similar comparison argument. Indeed, since
\[
-\,\mbox{div}\,\(|x|^{-\beta}\,\nabla w\)+|x|^{-\gamma}\,w^p\ge0\quad\forall\,x\in\R^d\setminus\{0\}
\]
and
\[
-\,\mbox{div}\,\(|x|^{-\beta}\,\nabla h\)+|x|^{-\gamma}\,h^p\le0\quad\forall\,x\in\R^d\setminus\{0\}
\]
if we choose $C$ such that $C^{p-1}\le(2-\gamma+\beta)\,\frac{d-\gamma-p\,(d-2-\beta)}{(p-1)^2}$, we easily see that
\[
w(x)\ge \( \min_{|x|=1}w(x) \wedge C \) |x|^{(\gamma-2-\beta)/(p-1)} \quad \forall x \in \mathbb{R}^d \setminus {B}_1 \, .
\]
This concludes the proof.
\end{proof}

\medskip Our next goal is to obtain growth and decay estimates, respectively, on the functions $\mathsf P$ and $u$ as they appear in the proof of Theorem~\ref{Thm:Rigidity} in Section~\ref{Sec:BE}, in order to prove Proposition~\ref{Prop:b}. We also need to estimate their derivatives near the origin and at infinity. Let us start by reminding the change of variables~\eqref{wv}, which in particular, by Lemma~\ref{Lem:estimates2}, implies that for some positive constants $C_1$ and $C_2$,
\[\label{asympv}
C_1\,s^{2/(1-p)}\le v(s,\omega)\le C_2\,s^{2/(1-p)}\quad\mbox{as}\; s\to+\infty\,.
\]
Then we perform the Emden-Fowler transformation
\be{Emden-Fowler}
v(s,\omega)=s^a\,\varphi(z,\omega)\quad\mbox{with}\quad z=-\log s\,,\quad a=\frac{2-n}2\,,
\ee
and see that $\varphi$ satisfies the equation
\be{varphieq}
-\,\alpha^2\,\varphi''-\Delta_\omega\,\varphi+a^2\alpha^2\varphi = e^{((n-2)\,p-n)\,z}\,\varphi^{2p-1} - e^{((n-2)\,p-n-2)\,z/2}\,\varphi^p=:h\quad\mbox{in}\;\mathcal C:=\R\times\S^{d-1}\ni(z,\omega)\,.
\ee
{}From here on we shall denote by $'$ the derivative with respect either to $z$ or to $s$, depending on the argument. By definition of $\varphi$ and using Lemma~\ref{Lem:estimates2}, we obtain that
\[\label{asympvarphi}
\varphi(z,\omega)\sim e^{\big(\frac{2-n}2+\frac2{p-1}\big)\,z}\quad\mbox{as}\; z\to -\infty\,,
\]
where we say that $f(z,\omega)\sim g(z,\omega)$ as $z\to+\infty$ (resp.~$z\to-\infty$) if the ratio $f/g$ is bounded from above and from below by positive constants, independently of $\omega$, and for $z$ (resp.~$-z$) large enough. Concerning $z\to+\infty$, we first note that Lemma~\ref{Lem:estimates1} and~\eqref{Emden-Fowler} show that $\varphi(z,\omega)\le O(e^{a\,z})$. The lower bound can be established by a comparison argument as in~\cite[Proposition~A.1]{DEL2015}, after noticing that $|h(z,\omega)|\le O(e^{(a-2)z})$. Hence we obtain that
\[
\varphi(z,\omega)\sim e^{a\,z}=e^{\frac{2-n}2\,z}\quad\mbox{as}\; z\to+\infty\,.
\]
Moreover, uniformly in $\omega$, we have that
\[
|h(z,\omega)| \le O\big( e^{-\frac{n+2}2\,z} \big) \quad\mbox{as}\; z\to+\infty\,,\quad |h(z,\omega)|\sim e^{\(-\frac{n+2}2+\frac{2\,p}{p-1}\)z}\quad\mbox{as}\; z\to -\infty\,,
\]
which in particular implies
\[\label{littleo}
|h(z,\omega)|=o\big(\varphi(z,\omega)\big)\quad\mbox{as}\; z\to+\infty\quad\mbox{and}\quad |h(z,\omega)|\sim\varphi(z,\omega)\quad\mbox{as}\; z\to-\infty\,.
\]
Finally, using~\cite[Theorem~8.32, p.~210]{MR1814364} on local $C^{1,\delta}$ estimates, as $|z|\to+\infty$ we see that all first derivatives of $\varphi$ converge to~$0$ at least with the same rate as $\varphi$. Next,~\cite[Theorem~8.10, p.~186]{MR1814364} provides local $\mathrm W^{k+2,2}$ estimates which, together with~\cite[Corollary~7.11, Theorem~8.10, and Corollary~8.11]{MR1814364}, up to choosing $k$ large enough, prove that
\be{asympt-beh-ders2}
|\varphi'(z,\omega)|\kern1.2pt,\;|\varphi''(z,\omega)|\kern1.2pt,\;|\nabla_\omega\varphi(z,\omega)|\kern1.2pt,\;|\nabla_\omega\varphi'(z,\omega)|\kern1.2pt,\;|\nabla_\omega\varphi''(z,\omega)|\kern1.2pt,\;|\Delta_\omega\,\varphi(z,\omega)|\le O(\varphi(z,\omega))\,,
\ee
uniformly in $\omega$. Here we denote by $\nabla_\omega$ the differentiation with respect to $\omega$. As a consequence, we have, uniformly in $\omega$, and for $\ell\in\{0,1,2\},\; t\in \{0,1\}$,
\be{asymph}
|\partial_z^\ell \nabla_\omega^t h(z,\omega)| \le O\big( e^{-\frac{n+2}2\,z} \big) \quad\mbox{as}\; z\to+\infty\,,\quad |\partial_z^\ell \nabla_\omega^t h(z,\omega)|\le O( e^{\(-\frac{n+2}2+\frac{2\,p}{p-1}\)z})\quad\mbox{as}\; z\to -\infty\,,
\ee
\be{asymphother}
|\Delta_\omega h(z,\omega)| \le O\big( e^{-\frac{n+2}2\,z} \big) \quad\mbox{as}\; z\to+\infty\,,\quad |\Delta_\omega h(z,\omega)|\le O( e^{\(-\frac{n+2}2+\frac{2\,p}{p-1}\)z})\quad\mbox{as}\; z\to -\infty\,.
\ee
\par\smallskip\begin{lemma}\label{Lem:decayinRd}{\sl Let $\beta$, $\gamma$ and $p$ satisfy the relations~\eqref{parameters} and assume $\alpha\le\alpha_{\rm FS}$. For any positive solution $w$ of~\eqref{ELeq} satisfying~\eqref{energyintegrals}, the \emph{pressure function} $\mathsf P=\frac m{1-m}\,u^{m-1}$ is such that $\mathsf P''$, $\mathsf P'/s$, $\mathsf P/s^2$, $\nabla_\omega\mathsf P'/s$, $\nabla_\omega\mathsf P/s^2$ and $\L\mathsf P$ are of class $C^\infty$ and bounded as $s\to+\infty$. On the other hand, as $s\to0_+$ we have
\begin{enumerate}
\item[(i)] $\isph{|\mathsf P'(s,\omega)|^2}\le O(1)$,
\item[(ii)] $\isph{|\nabla_\omega\mathsf P(s,\omega)|^2}\le O(s^2)$,
\item[(iii)] $\isph{|\mathsf P''(s,\omega)|^2}\le O(1/s^2)$,
\item[(iv)] $\isph{\left|\nabla_\omega\mathsf P'(s,\omega)-\tfrac1s\,\nabla_\omega\mathsf P(s,\omega)\right|^2}\le O(1)$,
\item[(v)] $\isph{\left| \frac{1}{s^2}\,\Delta_\omega \mathsf P(s,\omega)\right|^2}\le O(1/s^2)$.
\end{enumerate}}
\end{lemma}\par\smallskip
\begin{proof} By using the change of variables~\eqref{Emden-Fowler}, we see that
\[\label{Pvarphi}
\mathsf P(s,\omega)= \tfrac{p+1}{p-1}\,e^{-\frac12\,(n-2)\,(p-1)\,z}\,\varphi^{1-p}(z,\omega)\,,\quad z=-\log s\,.
\]
{}From~\eqref{asympt-beh-ders2} we easily deduce that uniformly in $\omega$, $\mathsf P''$, $\mathsf P'/s$, $\mathsf P/s^2$, $\nabla_\omega\mathsf P'/s$, $\nabla_\omega\mathsf P/s^2$ and $\L\mathsf P$ are of class $C^\infty$ and bounded as $s\to+\infty$. Moreover, as $s\to0_+$, we obtain that
\[
\big|\mathsf P'(s,\omega)\big|\le O\(\frac1s\(\frac{\varphi'(z,\omega)}{\varphi(z,\omega)}-a\)\)\quad\mbox{and}\quad\Big|\frac1s\,\nabla_\omega\mathsf P(s,\omega)\Big|\le O\(\frac1s\,\(\frac{\nabla_\omega\varphi(z,\omega)}{\varphi(z,\omega)}\)\)
\]
are of order at most $1/s$ uniformly in $\omega$. Similarly we obtain that
\begin{eqnarray*}
&&|\mathsf P''(s,\omega)|\le O\(\frac1{s^2}\(\frac{\varphi''(z,\omega)}{\varphi(z,\omega)}-\,p\,\frac{|\varphi'(z,\omega)|^2}{|\varphi(z,\omega)|^2}+\,\big(1-2\,a\,(1-p)\big)\,\frac{\varphi'(z,\omega)}{\varphi(z,\omega)}+a^2\,(1-p)-a\)\)\,,\\
&&\left|\frac{\nabla_\omega\mathsf P'(s,\omega)}{s}-\frac{a(1-p)}{s^2}\,\nabla_\omega\mathsf P(s,\omega)\right|\le O\(\frac1{s^2}\(\frac{\nabla_\omega\varphi'(z,\omega)}{\varphi(z,\omega)}-\frac{p\,\varphi'(z,\omega)\,\nabla_\omega\varphi(z,\omega)}{|\varphi(z,\omega)|^2}\)\)\,,\\
&&\frac1{s^2}\,|\Delta_\omega\,\mathsf P(s,\omega)|\le O\(\frac1{s^2}\,\(\frac{\Delta_\omega\,\varphi(z,\omega)}{\varphi(z,\omega)}-\,p\,\frac{|\nabla_\omega\varphi(z,\omega)|^2}{|\varphi(z,\omega)|^2}\)\)\,,
\end{eqnarray*}
are at most of order $1/s^2$ uniformly in $\omega$. This shows that $|\mathsf b(s)|\le O(s^{n-4})$ as $s\to0_+$ and concludes the proof if $4\le d<n$. When $d=2$ or $3$ and $n \le 4$, more detailed estimates are needed. Properties (i)--(v) amount to prove that
\begin{enumerate}
\item[(i)] $\isph{\left|\frac{\varphi'(z,\omega)}{\varphi(z,\omega)}-a\right|^2}\le O(e^{-2\,z})$,
\item[(ii)] $\isph{\left|\frac{\nabla_\omega\varphi(z,\omega)}{\varphi(z,\omega)}\right|^2}\le O(e^{-2\,z})$,
\item[(iii)] $\isph{\left|\frac{\varphi''(z,\omega)}{\varphi(z,\omega)}-\,p\,\frac{|\varphi'(z,\omega)|^2}{|\varphi(z,\omega)|^2}+\,\big(1-2\,a\,(1-p)\big)\frac{\varphi'(z,\omega)}{\varphi(z,\omega)}+a^2\,(1-p)-a\right|^2}\le O(e^{-2\,z})$,
\item[(iv)] $\isph{\left|\frac{\nabla_\omega\varphi'(z,\omega)}{\varphi(z,\omega)}-\frac{p\,\varphi'(z,\omega)\,\nabla_\omega\varphi(z,\omega)}{|\varphi(z,\omega)|^2}\right|^2}\le O(e^{-2\,z})$,
\item[(v)] $\isph{\left|\frac{\Delta_\omega\,\varphi(z,\omega)}{\varphi(z,\omega)}-\,p\,\frac{|\nabla_\omega\varphi(z,\omega)|^2}{|\varphi(z,\omega)|^2}\right|^2}\le O(e^{-2\,z})$,
\end{enumerate}
as $z\to+\infty$.

\medskip\noindent\emph{Step 1: Proof of\/ (ii) and (iv)}. If $w$ is a positive solution of~\eqref{ELeq}, then $\varphi$ is a positive solution to~\eqref{varphieq}. 
With $\ell\in\{0,1,2\}$, applying the operator $\nabla_\omega\partial_z^\ell$ to the equation~\eqref{varphieq} we obtain
\[
-\,\alpha^2\,(\nabla_\omega\partial_z^\ell\varphi)''-\,\nabla_\omega\,\Delta_\omega\,\partial_z^\ell\varphi+a^2\,\alpha^2\,\nabla_\omega\partial_z^\ell\varphi=\nabla_\omega\partial_z^\ell h(z,\omega)\quad\mbox{in}\;\mathcal C\,.
\]
Define 
\[
\chi_\ell(z):=\frac12\isph{|\nabla_\omega\partial_z^\ell\varphi|^2}\,,
\]
which by~\eqref{asympt-beh-ders2} converges to $0$ as $z\to \pm\infty$. Assume first that $\chi_\ell$ is a positive function.
After multiplying the above equation by $\nabla_\omega\partial_z^\ell\varphi$, integrating over $ \mathbb{S}^{d-1}$, integrating by parts and using
\[
\chi_\ell'=\isph{\nabla_\omega\partial_z^\ell\varphi\,\nabla_\omega\partial_z^\ell\varphi'}
\]
and
\[
\chi_\ell''=\isph{\nabla_\omega\partial_z^\ell\varphi\,\nabla_\omega\partial_z^\ell\varphi''}+\isph{|\nabla_\omega\partial_z^\ell\varphi'|^2}\,,
\]
we see that $\chi_\ell$ satisfies 
\[
-\,\chi_\ell''+\isph{|\nabla_\omega\partial_z^\ell\varphi'|^2} +\frac1{\alpha^2}\(\isph{|\Delta_\omega\,\partial_z^\ell\varphi|^2}-\lambda_1\isph{|\nabla_\omega\partial_z^\ell\varphi|^2}\)+ 2\,\(a^2+\frac{\lambda_1}{\alpha^2}\)\,\chi_\ell=\frac{h_\ell}{\alpha^2}\,,
\]
with $h_\ell :=\isph{\nabla_\omega\partial_z^\ell h\,\nabla_\omega\partial_z^\ell\varphi}$.
Then, using $\isph{\nabla_\omega\partial_z^\ell\varphi}=0$, by the Poincar\'e inequality we deduce
\[
\isph{|\Delta_\omega\,\partial_z^\ell\varphi|^2}\ge\lambda_1\isph{|\nabla_\omega\partial_z^\ell\varphi|^2}
\]
as \emph{e.g.}~in~\cite[Lemma~7]{MR3229793}, where $ \lambda_1:=d-1 $. A Cauchy-Schwarz inequality implies that
\[
-\,\chi_\ell''+\frac{|\chi_\ell'|^2}{2\,\chi_\ell}+\,2\,\(a^2+\frac{\lambda_1}{\alpha^2}\)\,\chi_\ell\le\frac{|h_\ell|}{\alpha^2} \, .
\]
The function $\zeta_\ell:=\sqrt{\chi_\ell}$ satisfies
\[
-\,\zeta_\ell''+\,\(a^2+\frac{\lambda_1}{\alpha^2}\)\,\zeta_\ell\le\frac{|h_\ell|}{2\,\alpha^2\,\zeta_\ell}\,.
\]
By the Cauchy-Schwarz inequality and~\eqref{asymph} we infer that $|h_\ell/\zeta_\ell| = O\big(e^{(a-2)\,z}\big)$ for $z\to+\infty$, and $|h_\ell/\zeta_\ell| = O\big(e^{(a+2/(p-1))\,z}\big)$ for $z\to-\infty$. By a simple comparison argument based on the Maximum Principle, and using the convergence of $\chi_\ell$ to $0$ at $\pm\infty$, we infer that
\[
\zeta_\ell(z)\le -\frac{e^{-\nu\,z}}{2\,\nu\,\alpha^2}\int_{-\infty}^ze^{\nu\,t}\,\frac{|h_\ell(t)|}{\zeta_\ell(z)}\,\dt-\frac{e^{\nu\,z}}{2\,\nu\,\alpha^2}\int_z^\infty e^{-\nu\,t}\,\frac{|h_\ell(t)|}{\zeta_\ell(z)}\,\dt
\]
if $\nu:=\sqrt{a^2+\lambda_1/{\alpha^2}}$. This is enough to deduce that $\zeta_\ell(z)\le O\big(e^{(a-1) z}\big)$ as $z\to+\infty$ after observing that the condition
\[\label{luckycomparison}
-\nu = -\sqrt{a^2+\lambda_1/{\alpha^2}} \le a-1
\]
is equivalent to the inequality $\alpha\le\alpha_{\rm FS}$. Hence we have shown that if $\chi_\ell$ is a positive function, then for $\alpha\le\alpha_{\rm FS}$,
\be{lestimates}
\chi_\ell(z)\le O\big(e^{\kern 0.5pt 2\,(a-1)\,z}\big)\quad\mbox{as}\; z\to+\infty\,.
\ee
In the case where $\chi_\ell$ is equal to $0$ at some points of $\R$, it is enough to do the above comparison argument on maximal positivity intervals of $\chi_\ell$ to deduce the same asymptotic estimate. Finally we observe that $\varphi(z,\omega)\sim e^{a\,z}$ as $z\to+\infty$, which ends the proof of (ii) considering the above estimate for $\chi_\ell$ when $\ell=0$. Moreover, the same estimate for $\ell=1$ together with (ii) and~\eqref{asympt-beh-ders2} proves (iv).

\medskip\noindent\emph{Step 2: Proof of\/(v)}. By applying the operator $\Delta_\omega$ to~\eqref{varphieq}, we obtain
\[
-\,\alpha^2\,(\Delta_\omega\,\varphi)''-\,\Delta_\omega^2\,\varphi+a^2\,\alpha^2\,\Delta_\omega\,\varphi=\Delta_\omega\,h \quad\mbox{in}\;\mathcal C\,.
\]
We proceed as in Step 1. With similar notations, by defining
\[
\chi_5(z):=\frac12\isph{|\Delta_\omega\,\varphi|^2}\,,
\]
after multiplying the equation by $\Delta_\omega\,\varphi$ and using the fact that
\[
-\isph{\Delta_\omega\,\varphi\,\Delta_\omega^2\,\varphi}=\isph{|\nabla_\omega\Delta_\omega\,\varphi|^2}\ge\lambda_1\isph{|\Delta_\omega\,\varphi|^2}\,,
\]
we obtain
\[
-\,\chi_5''+\frac{|\chi_5'|^2}{2\,\chi_\ell}+\,2\,\(a^2+\frac{\lambda_1}{\alpha^2}\)\,\chi_5\le\frac{|h_5|}{\alpha^2}
\]
with $h_5:=\isph{\Delta_\omega\,h\,\Delta_\omega\,\varphi}$. Again using the same arguments as above, together with~\eqref{asymphother}, we deduce that 
\[
\chi_5(z)\le O\big(e^{\kern 0.5pt 2\,(a-1)\,z}\big)\quad\mbox{as}\; z\to+\infty\,.
\]
This ends the proof of (v), using (ii),~\eqref{asympt-beh-ders2} and noticing again that $\varphi(z,\omega)\sim e^{a\,z}$ as $z\to+\infty$.

\medskip\noindent\emph{Step 3: Proof of (i) and (iii).}
Let us consider a positive solution $\varphi$ to~\eqref{varphieq} and define on $\R$ the function
\[
\varphi_0(z) := \frac{1}{\left| \mathbb{S}^{d-1} \right|} \isph{\varphi(z,\omega)}\,.
\]
By integrating~\eqref{varphieq} on $\S^{d-1}$, we know that $\varphi_0$ solves
\[
-\,\varphi_0''+a^2\,\varphi_0=\frac1{\alpha^2 \left| \mathbb{S}^{d-1} \right|}\,\isph{h(z,\omega)}=:\frac{h_0(z)}{\alpha^2 }\quad\forall\,z\in\R\,,
\]
with
\[
|h_0(z)| \le O \big( e^{-\frac{n+2}2\,z} \big) \quad\mbox{as}\; z\to+\infty\,,\quad |h_0(z)|\sim e^{\(-\frac{n+2}2+\frac{2\,p}{p-1}\)\,z}\quad\mbox{as}\; z\to -\infty\,.
\]
{}From the integral representation
\[
\varphi_0(z)=-\frac{e^{a\,z}}{2\,a\,\alpha^2}\int_{-\infty}^ze^{-at}\,h_0(t)\,\dt-\frac{e^{-a\,z}}{2\,a\,\alpha^2}\int_z^\infty e^{at}\,h_0(t)\,\dt\,,
\]
we deduce that as $z\to+\infty$, $\varphi_0(z) \sim e^{a\,z}$ and
\[
\frac{\varphi_0'(z)-a\,\varphi_0(z)}{\varphi(z,\omega)}\sim\,e^{-2a\,z}\int_z^{\infty}e^{at}\,h_0(t)\,\dt = O(e^{-2\,z})\,.
\]

If we define the function $\psi(z,\omega):=e^{-a\,z}\,\big(\varphi(z,\omega)-\varphi_0(z)\big)$, we may observe that it is bounded for $z$ positive and moreover
\[
\frac{\varphi'(z,\omega)}{\varphi(z,\omega)}-a=O(e^{-2\,z})+\frac{\psi'(z,\omega)}{e^{-a\,z}\,\varphi(z,\omega)}\quad\mbox{as}\;\; z\to\,+\infty\,.
\]
We recall that $e^{-a\,z}\,\varphi(z,\omega)$ is bounded away from $0$ by a positive constant as $z\to+\infty$. Hence we know that
\be{dzpsi}
\Big|\frac{\varphi'(z,\omega)}{\varphi(z,\omega)}-a\Big|\le O\(|\psi'(z,\omega)\)+O(e^{-2\,z})\,.
\ee
By the Poincar\'e inequality and estimate~\eqref{lestimates} with $\ell=0$, we have
\[
\isph{|\psi|^2}= e^{-2az} \isph{|\varphi-\varphi_0|^2}\le\frac{e^{-2az}}{\lambda_1} \isph{|\nabla_\omega\varphi|^2}\le O(e^{-2z})\,.
\]
Moreover, by the estimate~\eqref{lestimates} with $\ell=1$, we also obtain
\[
e^{-2az} \isph{|\varphi'-\varphi'_0|^2}\le\frac{e^{-2az}}{\lambda_1} \isph{|\nabla_\omega\varphi'|^2}\le O(e^{-2z})\,.
\]
Hence, since $\psi' = -\,a\,\psi+ e^{-az}\,(\varphi'-\varphi'_0)$, the above estimates imply that 
\[
\isph{|\psi|^2} \, + \, \isph{|\psi'|^2}\le O(e^{-2z})\,,
\]
which together with~\eqref{dzpsi} ends the proof of (i).

To prove (iii), we first check that
\[
\frac{\varphi''}{\varphi}-\,p\,\frac{|\varphi'|^2}{|\varphi|^2}+\,\big(1-2\,a\,(1-p)\big)\frac{\varphi'}{\varphi}+a^2\,(1-p)-a= O(|\psi'|+|\psi'|^2+|\psi''|)+O(e^{-2\,z})\,,
\]
and so it remains to prove that $\isph{|\psi''|^2}$ is of order $O(e^{-2\,z})$. Since 
\[
\psi'' = a^2\,\psi -\,2\,a\,e^{-az}\,(\varphi'-\varphi'_0) + e^{-az}\,(\varphi''-\varphi''_0) \,,
\]
using the above estimates, we have only to estimate the term with the second derivatives. This can be done as above by the Poincar\'e inequality,
\[
e^{-2az} \isph{|\varphi''-\varphi''_0|^2}\le\frac{e^{-2az}}{\lambda_1} \isph{|\nabla_\omega\varphi''|^2}\le O(e^{-2z})\,,
\]
based on the estimate~\eqref{lestimates} with $\ell=2$.
This ends the proof of (iii).
\end{proof}

\medskip\noindent\emph{Proof of Proposition~\ref{Prop:b}}
It is straightforward to verify that the boundedness of $\mathsf P''$, $\mathsf P'/s$, $\mathsf P/s^2$, $\nabla_\omega\mathsf P'/s$, $\nabla_\omega\mathsf P/s^2$, $\L\mathsf P$ as $ s \to +\infty $ and the integral estimates (i)-(v) as $ s \to 0^+ $ from Lemma~\ref{Lem:decayinRd} are enough in order to establish~\eqref{FisherDerBC},~\eqref{FisherBC} and~\eqref{MassBC}.
\hfill \ \qed 

\bigskip\noindent{\small{\bf Acknowlegments.} This research has been partially supported by the projects \emph{STAB} (J.D.) and \emph{Kibord} (J.D.) of the French National Research Agency (ANR), and by the NSF grant DMS-1301555 (M.L.). J.D.~thanks the University of Pavia for support. M.L.~thanks the Humboldt Foundation for support. M.M.~has been partially funded by the National Research Project ``Calculus of Variations'' (PRIN 2010-11, Italy).
\par\smallskip\noindent\copyright~2016 by the authors. This paper may be reproduced, in its entirety, for non-commercial purposes.}

\newpage

\begin{thebibliography}{10}

\bibitem{Bakry-Emery85}
{\sc D.~Bakry and M.~{\'E}mery}, {\em Diffusions hypercontractives}, in
  S\'eminaire de probabilit\'es, XIX, 1983/84, vol.~1123 of Lecture Notes in
  Math., Springer, Berlin, 1985, pp.~177--206.

\bibitem{MR3155209}
{\sc D.~Bakry, I.~Gentil, and M.~Ledoux}, {\em Analysis and geometry of
  {M}arkov diffusion operators}, vol.~348 of Grundlehren der Mathematischen
  Wissenschaften [Fundamental Principles of Mathematical Sciences], Springer,
  Cham, 2014.

\bibitem{MR0188004}
{\sc N.~M. Blachman}, {\em The convolution inequality for entropy powers}, IEEE
  Trans. Information Theory, IT-11 (1965), pp.~267--271.

\bibitem{2016arXiv160208319B}
{\sc M.~{Bonforte}, J.~{Dolbeault}, M.~{Muratori}, and B.~{Nazaret}}, {\em
  {Weighted fast diffusion equations (Part I): Sharp asymptotic rates without
  symmetry and symmetry breaking in Caffarelli-Kohn-Nirenberg inequalities}}.
\newblock Preprint hal-01279326 \& ar{X}iv: 1602.08319, Feb. 2016.

\bibitem{Caffarelli-Kohn-Nirenberg-84}
{\sc L.~Caffarelli, R.~Kohn, and L.~Nirenberg}, {\em First order interpolation
  inequalities with weights}, Compositio {M}ath., 53 (1984), pp.~259--275.

\bibitem{MR1853037}
{\sc J.~A. Carrillo, A.~J{\"u}ngel, P.~A. Markowich, G.~Toscani, and
  A.~Unterreiter}, {\em Entropy dissipation methods for degenerate parabolic
  problems and generalized {S}obolev inequalities}, Monatsh. {M}ath., 133
  (2001), pp.~1--82.

\bibitem{MR1777035}
{\sc J.~A. Carrillo and G.~Toscani}, {\em Asymptotic {$L^1$}-decay of solutions
  of the porous medium equation to self-similarity}, Indiana Univ. Math. J., 49
  (2000), pp.~113--142.

\bibitem{MR1986060}
{\sc J.~A. Carrillo and J.~L. V{\'a}zquez}, {\em Fine asymptotics for fast
  diffusion equations}, Comm. Partial Differential Equations, 28 (2003),
  pp.~1023--1056.

\bibitem{Catrina-Wang-01}
{\sc F.~Catrina and Z.-Q. Wang}, {\em On the {C}affarelli-{K}ohn-{N}irenberg
  inequalities: sharp constants, existence (and nonexistence), and symmetry of
  extremal functions}, {C}omm. {P}ure {A}ppl. {M}ath., 54 (2001), pp.~229--258.

\bibitem{MR823597}
{\sc M.~H.~M. Costa}, {\em A new entropy power inequality}, IEEE Trans. Inform.
  Theory, 31 (1985), pp.~751--760.

\bibitem{MR1940370}
{\sc M.~Del~Pino and J.~Dolbeault}, {\em Best constants for
  {G}agliardo-{N}irenberg inequalities and applications to nonlinear
  diffusions}, {J}ournal de {M}ath\'ematiques {P}ures et {A}ppliqu\'ees. (9),
  81 (2002), pp.~847--875.

\bibitem{MR3229793}
{\sc J.~Dolbeault, M.~J. Esteban, and M.~Loss}, {\em Nonlinear flows and
  rigidity results on compact manifolds}, J. {F}unct. {A}nal., 267 (2014),
  pp.~1338--1363.

\bibitem{dolbeault:hal-01206975}
\leavevmode\vrule height 2pt depth -1.6pt width 23pt, {\em {Interpolation
  inequalities on the sphere: linear vs. nonlinear flows}}, {Annales de la
  Facult{\'e} des Sciences de Toulouse. Math{\'e}matiques},  (2016).

\bibitem{DEL2015}
\leavevmode\vrule height 2pt depth -1.6pt width 23pt, {\em Rigidity versus
  symmetry breaking via nonlinear flows on cylinders and {E}uclidean spaces},
  to appear in Inventiones {M}athematicae,
  http://link.springer.com/article/10.1007/s00222-016-0656-6 (2016).

\bibitem{dolbeault:hal-01286546}
\leavevmode\vrule height 2pt depth -1.6pt width 23pt, {\em {Symmetry of
  optimizers of the Caffarelli-Kohn-Nirenberg inequalities}}.
\newblock Preprint hal-01286546 \& ar{X}iv: 1603.03574, Mar. 2016.

\bibitem{0902}
{\sc J.~Dolbeault, M.~J. Esteban, M.~Loss, and G.~Tarantello}, {\em On the
  symmetry of extremals for the {C}affarelli-{K}ohn-{N}irenberg inequalities},
  Advanced Nonlinear Studies, 9 (2009), pp.~713--727.

\bibitem{DMN2015}
{\sc J.~Dolbeault, M.~Muratori, and B.~Nazaret}, {\em Weighted interpolation
  inequalities: a pertur\-bation approach}.
\newblock Preprint hal-01207009 \& ar{X}iv: 1509.09127, to appear in {M}ath.
  {A}nnalen, 2016.

\bibitem{1501}
{\sc J.~Dolbeault and G.~Toscani}, {\em Nonlinear diffusions: Extremal
  properties of {B}arenblatt profiles, best matching and delays}, Nonlinear
  Analysis: Theory, Methods \& Applications, 138 (2016), pp.~31 -- 43.
\newblock Nonlinear Partial Differential Equations, in honor of Juan Luis
  V{\'a}zquez for his 70th birthday.

\bibitem{Felli-Schneider-03}
{\sc V.~Felli and M.~Schneider}, {\em Perturbation results of critical elliptic
  equations of {C}affarelli-{K}ohn-{N}irenberg type}, {J}. {D}ifferential
  {E}quations, 191 (2003), pp.~121--142.

\bibitem{MR1814364}
{\sc D.~Gilbarg and N.~S. Trudinger}, {\em Elliptic partial differential
  equations of second order}, Classics in Mathematics, Springer-Verlag, Berlin,
  2001.
\newblock Reprint of the 1998 edition.

\bibitem{MR3200617}
{\sc G.~Savar{\'e} and G.~Toscani}, {\em The concavity of {R}\'enyi entropy
  power}, {IEEE} {T}rans. {I}nform. {T}heory, 60 (2014), pp.~2687--2693.

\bibitem{MR2282669}
{\sc J.~L. V{\'a}zquez}, {\em Smoothing and decay estimates for nonlinear
  diffusion equations. Equations of porous medium type}, vol.~33, Oxford
  University Press, Oxford, 2006.

\bibitem{MR1768665}
{\sc C.~Villani}, {\em A short proof of the ``concavity of entropy power''},
  IEEE Trans. Inform. Theory, 46 (2000), pp.~1695--1696.

\bibitem{MR479373}
{\sc F.~B. Weissler}, {\em Logarithmic {S}obolev inequalities for the
  heat-diffusion semigroup}, Trans. Amer. Math. Soc., 237 (1978), pp.~255--269.

\end{thebibliography}

\end{document}